%% file: Sobolev_arXiv_v2.tex
\documentclass[oneside,11pt]{article}

\input{preamble_arXiv.tex}
\newcommand{\pr}{^{\prime}}

\usepackage[pdftex,final,bookmarksnumbered,bookmarksopen=false,breaklinks,colorlinks]{hyperref}
\hypersetup{
	final,
	bookmarksnumbered=true,
	bookmarksopen=true,
	bookmarksopenlevel=0,
	unicode=false,          %
	pdftoolbar=true,        %
	pdfmenubar=true,        %
	pdffitwindow=false,     %
	pdftitle={On a class of Sobolev tests for symmetry of directions, their detection thresholds, and asymptotic powers},    %
	pdfdisplaydoctitle=true, %
	pdftoolbar=true,
	pdfmenubar=true,
	pdflang={English},
	pdfauthor={Eduardo Garcia-Portugues, Davy Paindaveine, Thomas Verdebout},     %
	pdfsubject={arXiv paper},   %
	pdfcreator={Eduardo Garcia-Portugues},   %
	pdfproducer={Eduardo Garcia-Portugues}, %
	pdfkeywords={Consistency rates}{Detection thresholds}{Directional statistics}{Hypothesis tests}, %
	pdfnewwindow=true,      %
	breaklinks=true,
	hidelinks,       
	linkcolor=black,          %
	citecolor=black,        %
	filecolor=magenta,      %
	urlcolor=cyan           %
}

\newif\ifmain
\maintrue
\newif\ifsupplement
\supplementtrue

\newif\iffigstabs
\figstabstrue

\begin{document}

\ifmain

%-----------------------------------------------%
\title{On a class of Sobolev tests for symmetry of directions, their detection thresholds, and asymptotic powers}
\setlength{\droptitle}{-1cm}
\predate{}%
\postdate{}%
\date{}
%-----------------------------------------------%

%-----------------------------------------------%
\author{Eduardo Garc\'ia-Portugu\'es$^{1}$, Davy Paindaveine$^{2,3}$, and Thomas Verdebout$^{2}$}
\footnotetext[1]{Department of Statistics, Universidad Carlos III de Madrid (Spain).}
\footnotetext[2]{D\'{e}partement de Math\'{e}matique and ECARES, Universit\'{e} libre de Bruxelles (Belgium).}
\footnotetext[3]{Corresponding author. e-mail: \href{mailto:dpaindav@ulb.ac.be}{dpaindav@ulb.ac.be}.}
\maketitle
%-----------------------------------------------%

\begin{abstract}
	We consider a class of symmetry hypothesis testing problems including testing isotropy on~$\mathbb{R}^d$ and testing rotational symmetry on the hypersphere~$\mathcal{S}^{d-1}$. For this class, we study the null and non-null behaviors of \emph{Sobolev tests}, with emphasis on their consistency rates. Our main results show that: (\textit{i})~Sobolev tests exhibit a \emph{detection threshold} \cite[see][]{Bha2019a,Bha2019} that does not only depend on the coefficients defining these tests; and (\textit{ii})~tests with non-zero coefficients at odd (respectively, even) ranks only are blind to alternatives with angular functions whose $k$th-order derivatives at zero vanish for any $k$~odd (even). Our non-standard asymptotic results are illustrated with Monte~Carlo exercises. A case study in astronomy applies the testing toolbox to evaluate the symmetry of orbits of long- and short-period comets.
\end{abstract}
\begin{flushleft}
\end{flushleft}

%-------------------------------%
\section{Introduction}
\label{sec:Intro}
%-------------------------------%

Directional statistics involve observations such as directions (unit vectors), axes (lines through the origin), rotations, etc., with values on manifolds---not, as in classical multivariate statistical analysis, unrestricted $\mathbb{R}^d$-valued vectors. Directional data can be found in a variety of fields, including astronomy \citep[e.g.,][]{Marinucci2008}, medicine \citep[e.g.,][]{Vuollo2016}, %
genetics \citep[e.g.,][]{Dortet-Bernadet2008}, %
to cite only a few. Directional data are often represented as points on the unit hypersphere~$ \mathcal{S}^{d-1}:= \{\xb \in \mathbb{R}^d: \|\xb \|^2 := \xb\pr \xb =1\}$ of $\mathbb{R}^d$. The survey in \cite{Pewsey2021} reviews recent advances in directional statistics.

An important problem within directional statistics is the problem of testing uniformity on the unit hypersphere~$\mathcal{S}^{d-1}$. Testing isotropy is one of the oldest problems in multivariate statistics and can be traced back to \cite{Bernoulli1735}. It is a fundamental problem that is still much considered nowadays. To cite only a few recent works, uniformity tests were studied in the context of noisy data in \cite{Lacour2014} and \cite{Kim2016}. Tests in the high-dimensional context were considered in \cite{Cai2012}, \cite{Cai2013}, and \cite{Cutting2017}. Projection-based tests were proposed in \cite{Cuesta-Albertos2009} and \cite{Garcia-Portugues2020b}. \emph{Sobolev tests} for the problem have been proposed in \cite{Beran1968,Beran1969}, \cite{Gine1975}, \cite{Bogdan2002}, and \cite{Jupp2008}. We refer to \cite{Garcia-Portugues2020a} for a review of uniformity~tests.

When uniformity of directions is rejected, practitioners often turn to the most popular models in directional statistics, namely the rotationally symmetric models. Rotational symmetry has been extensively adopted as a core assumption for performing inference with directional data. A (far from exhaustive) list of references illustrating this is as follows: \cite{Rivest1989}, \cite{Ko1993}, and \cite{Chang2001} considered regression and M-estimation under rotationally symmetric assumptions; \cite{Larsen2002} considered von Mises--Fisher likelihood ratios; \cite{Tsai2007}, \cite{Ley2013}, and \cite{Paindaveine2015} proposed rank tests and estimators for the mode of a rotationally symmetric distribution; \cite{Ley2014} proposed a concept of quantiles for rotationally symmetric distributions; \cite{Paindaveine2017} considered inference for the direction of weak rotationally symmetric signals, whereas \cite{Paindaveine2020a} tackled high-dimensional hypothesis testing in the same framework. Tests for rotational symmetry have been proposed, e.g., in \cite{Schach1969}, \cite{Pewsey2002a}, and \cite{Meintanis2019} for the circular case ($d=2$), and in \cite{Jupp1983}, \cite{Ley2017b}, and \cite{Garcia-Portugues2020} for the general case.

In the present paper, we propose a class of Sobolev tests for a wide class of symmetry testing problems. Particular instances of the tests we propose are the classical Sobolev tests of uniformity studied for instance in \cite{Jupp2008} and the tests for rotational symmetry proposed in \cite{Garcia-Portugues2020}. The asymptotic behavior of our Sobolev tests under the null hypothesis directly follows from classical results. However, this is far to be the case regarding their non-null asymptotic behavior, and, as a matter of fact, the non-null behavior of Sobolev tests remains largely unexplored even in the important particular case of testing uniformity. In this work, we therefore extensively study the asymptotic behavior of our Sobolev tests under natural sequences of local alternatives defined through a so-called \emph{angular function} $f$ and some positive sequence of concentration parameters $\kappa_n$; the larger $\kappa_n$, the more severe the deviation from the null hypothesis (throughout, $n$ is the sample~size). Our results reveal that a Sobolev test characterized by a sequence $v=(v_k)$ (see Section \ref{sec:Sobol} for more details) typically shows non-trivial asymptotic powers against alternatives in~$\kappa_n\sim n^{-1/(2k_v)}$, where~$k_v$ is the rank of the first non-zero coefficient in the sequence~$v=(v_k)$. More precisely, this is the case if~$f^{\underline{k_v}}(0)\neq 0$ (henceforth, $f^{\underline{k}}$ stands for the $k$th derivative of $f$). If this derivative is equal to zero, then slower consistency rates are actually achieved. Therefore, consistency rates of Sobolev tests do not only depend on the particular test used (that is, on the sequence~$v$) but also on which derivatives~$f^{\underline{k}}(0)$ are equal to zero. In the case where~$f^{\underline{k_v}}(0)=0$, we determine the local alternatives (if any) that can be detected by a Sobolev test. Remarkably, parity of the ranks~$k$ for which~$v_k$ and~$f^{\underline{k}}(0)$ are non-zero plays a key role in the resulting consistency rates. A relevant practical consequence of this parity feature is that any finite Sobolev test for which both $v_1$ and $v_2$ are non-zero will exhibit the same consistency rate as the oracle Sobolev test based on $v_k = \delta_{k,\ell}$, where $\ell$ is the smallest positive integer for which~$f^{\underline{\ell}}(0)\neq 0$.%

The subtle nature of our asymptotic results prevents us from being more specific in this introduction. Let us still mention, however, that whenever we identify in the sequel the consistency rate of a Sobolev test, we will provide the exact expression of the corresponding limiting non-trivial powers. In addition, these explicit limiting powers are practically relevant to determine ``indicative powers'' that inform on the confidence of the symmetry test in rejecting the null hypothesis given the estimated rotationally symmetric alternative. Due to its generality and indicative powers, the new testing class enables a deep assessment of the types of symmetry present in the orbits of long- and short-period comets.

The paper is organized as follows. In Section~\ref{sec:testprob}, we introduce the class of symmetry testing problems we tackle and we define the proposed Sobolev tests. In Section~\ref{sec:Sobol}, we first introduce the mathematical objects (including spherical harmonics as well as Chebyshev and Gegenbauer polynomials) that are needed to study the asymptotic properties of our Sobolev tests. We then recall what is known about the asymptotic null behavior of these tests in classical problems. In Section~\ref{sec:moments}, we study how the first moments of spherical harmonics evaluated at observations behave in the vicinity of uniformity, which requires investigating the first moments of random Chebyshev and Gegenbauer polynomials. While the results of Section~\ref{sec:moments} are of independent interest, we exploit them in this paper to identify the consistency rates of Sobolev tests and to derive the corresponding local asymptotic powers; this is done in Section~\ref{sec:localpowfiniteinf} for finite and infinite Sobolev tests (we say that a Sobolev test is finite if and only if finitely many of its coefficients~$v_k$ are non-zero). In Section~\ref{sec:simus}, we conduct several Monte Carlo exercises in order to illustrate the non-null results obtained in the previous sections. The usefulness of Sobolev tests is illustrated in Section \ref{sec:data} with a case study on the symmetry of comet orbits. Finally, we briefly discuss perspectives for future research in Section~\ref{sec:perspectives}. All proofs are provided in the Supplementary Material.

For convenience, we collect here some notation that will be used in the sequel. We will write~$\mathbb{I}[A]$ for the indicator function of the condition~$A$ (taking value one if~$A$ is satisfied and value zero otherwise) and we will denote as~${\bf e}_{\ell,k}$ the $\ell$th vector of the canonical basis of~$\R^k$. To make the notation lighter, we will simply write~${\bf e}_{\ell}$ whenever the dimension of the ambient Euclidean space can be recovered from the mathematical objects involved in the various formulae. We will denote as~$\nu_d$ the uniform probability measure on~$\mathcal{S}^{d-1}$, that is, $\sigma_d/\omega_{d-1}$, where~$\sigma_d$ is the surface area measure on~$\mathcal{S}^{d-1}$ and~$\omega_{d-1}=\sigma_d(\mathcal{S}^{d-1})=2\pi^{d/2}/\Gamma(d/2)$ is the surface area of~$\mathcal{S}^{d-1}$. This allows us to consider the inner product~$\langle f, g \rangle:=\int_{\mathcal{S}^{d-1}}f(\mathbf{x})g(\mathbf{x})\,\mathrm{d}\nu_d(\mathbf{x})$ on $L^2(\mathcal{S}^{d-1},\nu_d)$. For a real number~$a$ and a nonnegative integer~$k$, we will use the Pochhammer symbol $(a)_k:=a(a+1)\cdots(a+k-1)$ for~$k>0$ and~$(a)_k:=1$ for~$k=0$. The relation~$k\sim \ell$ on~$\mathbb{N}$ is satisfied when~$k$ and~$\ell$ share the same parity. We will write ${\rm diag}({\bf B}_1, \ldots, {\bf B}_\ell)$ for the block-diagonal matrix with diagonal blocks ${\bf B}_1, \ldots, {\bf B}_\ell$. Weak convergence will be denoted as~$\stackrel{\mathcal{D}}{\to}$ and is considered throughout the paper as the sample size~$n$ diverges to infinity. To simplify the exposition in the limiting distributions, we will recurrently write
$
\sum_{k=1}^\infty w_k \chi_{d_k}^2(\xi_k)
$
for denoting weighted sums of mutually independent random variables that have chi-squared distribution with degrees of freedom $d_k$ and (possible) non-centrality parameters $\xi_k$ (with $\chi_{d_k}^2\equiv\chi_{d_k}^2(0)$), for any $k\geq 1$. Finally, the notation
\begin{align}
	\label{cpormula}
	c_d
	:=
	\frac{\Gamma(\frac{d}{2})}{\sqrt{\pi}\Gamma(\frac{d-1}{2})}= 1 \, \Big/\int_{-1}^1 (1-s^2)^{(d-3)/2} \,\mathrm{d}s
\end{align}
will be used throughout the paper.

%-------------------------------%
\section{Testing problems and the proposed class of Sobolev tests}
\label{sec:testprob}
%-------------------------------%

Assume that we observe $n$ independent copies $\Zb_1, \ldots, \Zb_n$ of a $d$-dimensional random vector~$\Zb$ and let $\Gamb$ be a $p \times d$ semi-orthogonal matrix in the sense that
$$
\Gamb\Gamb\pr = {\bf I}_p \quad {\rm and} \quad \Gamb\pr \Gamb= \sum_{j=1}^p {\pmb \mu}_j {\pmb \mu}_j\pr,
$$
for $d$-dimensional orthonormal vectors ${\pmb \mu}_1, \ldots, {\pmb \mu}_p$ (throughout, $\mathbf{I}_k$ denotes the $k$-dimensional identity matrix). In this work, we consider the problem of testing the null hypothesis
\begin{align}
	{\cal H}_{0}: \frac{\Gamb \Zb}{\| \Gamb \Zb \|} \sim {\rm Unif}({\cal S}^{p-1}),\label{eq:H0}
\end{align}
where ${\rm Unif}({\cal S}^{p-1})$ stands for the uniform probability measure over~$\mathcal{S}^{p-1}$. Letting $\Ub_i (\Gamb):= \Gamb \Zb_i / \| \Gamb \Zb_i \|$, $i=1, \ldots, n$, we have two especially relevant instances for $\Gamb$ and thus \eqref{eq:H0}:
\begin{enumerate}[label=(\textit{\roman{*}}),ref=\textit{\roman{*}}]
	\item taking $\Gamb= {\bf I}_d$ ($p=d$) yields the classical problem of testing the null hypothesis that~$\Zb_1 / \| \Zb_1 \|,\allowbreak \ldots,\Zb_n / \| \Zb_n \|$ form a random sample from the uniform distribution on $\mathcal{S}^{d-1}$ (equivalently, that the~$\Zb_i$'s are sampled from an isotropic distribution on $\mathbb{R}^d$);
	\item taking $\Gamb$ such that $\Gamb\pr \Gamb= {\bf I}_d- {\pmb \theta} {\pmb \theta}\pr$ for some ${\pmb \theta} \in \mathcal{S}^{d-1}$ ($p=d-1$) yields the problem of testing that the $\Ub_i (\Gamb)$'s are uniformly distributed over $\mathcal{S}^{d-2}$, which allows one to test the null hypothesis that $\Zb_1, \ldots, \Zb_n$ share a common rotationally symmetric distribution about ${\pmb \theta}$; if the latter null hypothesis holds, then the $\Ub_i (\Gamb)$'s are indeed uniformly distributed over $\mathcal{S}^{d-2}$; this is the rationale used in the construction of the tests for rotational symmetry proposed in \cite{Garcia-Portugues2020}.
\end{enumerate}

To test the uniformity of the $\Ub_i (\Gamb)$'s, a very natural idea is to consider a \cite{Rayleigh1919}-type test, which rejects the null hypothesis at asymptotic level~$\alpha$ when
\begin{align}
	\label{Raytest}
	np \| \bar{\Ub} (\Gamb) \|^2 > \chi_{p;1-\alpha}^2,
\end{align}
where $\bar{\Ub} (\Gamb):= n^{-1}\sum_{i=1}^n \Ub_i (\Gamb)$ is the sample average of the observations and $\chi_{\ell;\beta}^2$ denotes the order-$\beta$
quantile of the chi-square distribution with~$\ell$ degrees of freedom. As explained in \citet[Section 6.3.1]{Mardia1999a}, the classical Rayleigh test for uniformity over $\mathcal{S}^{p-1}$ is the likelihood ratio test for testing uniformity within the rotationally symmetric von Mises--Fisher model characterized by densities on~$\mathcal{S}^{p-1}$ (throughout, densities are with respect to the surface area measure on~$\mathcal{S}^{p-1}$) proportional to~$
\ub \mapsto %
\exp(\kappa \ub\pr \thetab)$; here,~$\thetab \in {\cal S}^{p-1}$ is a location parameter and $\kappa \geq 0$ is a concentration parameter (the larger~$\kappa$, the more the probability mass is concentrated about~$\thetab$).
Clearly, the Rayleigh-type test in~\eqref{Raytest} is designed to show power under alternatives for which~${\rm E}[\Ub_1(\Gamb)]$ deviates from its null value (namely, the zero $p$-vector), but it will be blind to other alternatives. In particular, it will be blind to \mbox{(centro-)}symmetric (i.e., ``axial'') alternatives. To detect such symmetric alternatives, first-order moments can be replaced with second-order ones, which leads to the \cite{Bingham1974} test rejecting the null hypothesis of uniformity at asymptotic level $\alpha$ when
\begin{align}
	\label{Bing}
	\frac{np (p+2)}{2}
	{\rm tr}\bigg[ \bigg(\Sb (\Gamb) - \frac{1}{p}\mathbf{I}_p\bigg)^2 \bigg]
	>
	\chi_{(p-1)(p+2)/2;1-\alpha}^2
	,
\end{align}
where $\Sb (\Gamb):= n^{-1} \sum_{i=1}^{n} \Ub_i (\Gamb) \Ub_i\pr (\Gamb)$ is the empirical covariance matrix of the observations, relative to the fixed location vector zero.

The Rayleigh- and Bingham-type tests in \eqref{Raytest} and \eqref{Bing} have already been extensively studied in the literature:
\begin{enumerate}[label=(\textit{\roman{*}}),ref=\textit{\roman{*}}]
	\item for $\Gamb= {\bf I}_d$ (the classical uniformity problem), \cite{Chikuse2003} showed that the Rayleigh test is a locally most powerful rotation-invariant test against von Mises--Fisher alternatives. \cite{Jupp2001} proposed modifications of the Rayleigh and Bingham tests to improve their asymptotic chi-square approximations. \cite{Cutting2017}
	\vspace{-.4mm}
	established that the Rayleigh test detects local von Mises--Fisher alternatives in $\kappa_n \sim p_n^{3/4}/ \sqrt{n}$ in the high-dimensional case. The Bingham test is the likelihood ratio test for testing uniformity against the alternatives associated with the \cite{Bingham1974} distributions or the angular Gaussian distributions from \cite{Tyler1987}. \cite{Cutting2020} derived the local powers of the Bingham test under rotationally symmetric alternatives. The consistency rate of this test in high dimensions was also recently obtained in \cite{Cutting2021}. Finally, note that the Rayleigh and Bingham tests can be seen as the spatial sign test for location and for sphericity, respectively; see \cite{HP02} and \cite{HP06}, respectively.
	\item For $\Gamb$ such that $\Gamb\pr \Gamb= {\bf I}_d- {\pmb \mu}_1 {\pmb \mu_1}\pr$,   ${\pmb \mu}_1 \in \mathcal{S}^{d-1}$, the tests for rotational symmetry proposed in \cite{Garcia-Portugues2020} are the Rayleigh and Bingham tests in \eqref{Raytest} and \eqref{Bing}, respectively.  \label{RayBingRot}
\end{enumerate}

In this work, we propose Sobolev tests that reject %
\eqref{eq:H0} for large values~of %
\begin{align}
	\label{Sobol1}
	S^{(n)}_v
	:=
	\frac{1}{n} \sum_{i,j=1}^n\sum_{k=1}^\infty v_k^2
	h_{p,k}(\Ub_i'(\Gamb)\Ub_j (\Gamb))
	,
\end{align}
where
$
h_{p,k}(s)
:=
(
2\mathbb{I}[p=2]
+
\big(1+ \frac{2k}{p-2}\big) \mathbb{I}[p\geq 3]
)
C_k^{(p-2)/2}(s)
$, $s\in[-1,1]$, involves Chebyshev polynomials of the first kind~$C_k^{0}$ for~$p=2$ and Gegenbauer polynomials~$C_k^{\lambda}$, $\lambda>0$, for~$p\geq 3$; see Section~\ref{sec:Sobol} for a definition of these polynomials, for the construction of Sobolev tests, and for the summability condition the sequence of coefficients~$v=(v_k)$ should obviously satisfy. The Rayleigh and Bingham-type tests are obtained by taking $v_k=\delta_{k,1}$ and~$v_k=\delta_{k,2}$, where~$\delta_{k,\ell}$ is the usual Kronecker delta and therefore belong to the class of Sobolev tests in~\eqref{Sobol1}. The class of tests we propose therefore (\textit{i}) contains all the classical Sobolev tests of uniformity on $\mathcal{S}^{d-1}$ obtained with $p=d$ and (\textit{ii}) contains with $p=d-1$ (but is not limited to) the Rayleigh- and Bingham-type tests of rotational symmetry studied in \cite{Garcia-Portugues2020}. Our tests can also be used to perform tests for more complicated symmetry structures with $p < d-1$: in that case, they can be seen as subsphericity~tests.

In the next section, we summarize the asymptotic behavior of our test statistics under the null hypothesis and then we tackle their nonstandard asymptotic behavior under local alternatives in the sections that follow.

%

%

%

%
%-------------------------------%
\section{Spherical harmonics, Sobolev tests, and behavior under the null hypothesis }
\label{sec:Sobol}
%-------------------------------%

Sobolev tests on $\mathcal{S}^{p-1}$ are constructed via an orthonormal basis of $L^2(\mathcal{S}^{p-1},\nu_p)$ associated with \emph{spherical harmonics}. In this section, we first briefly present spherical harmonics, then we explain how they lead to Sobolev tests, and finally we review what is known on the asymptotic null behavior of these tests.

Fix an integer~$p\geq 2$ and a nonnegative integer~$k$. Denote as~$\mathcal{P}^p_k$ the space of real homogeneous polynomials of degree~$k$ and as~$\mathcal{H}^p_k$ the linear subspace of~$\mathcal{P}^p_k$ collecting polynomials whose Laplacian is zero.
Spherical harmonics of degree~$k$ simply are the restrictions to~$\mathcal{S}^{p-1}$ of polynomials in~$\mathcal{H}^p_k$. In the sequel, we will also denote as~$\mathcal{H}^p_k$ the collection of spherical harmonics of degree~$k$. It is well-known that~$\mathcal{H}^p_k$ is a vector space with dimension
$$
d_{p,k}
=
\binom{p+k-1}{k}-\binom{p+k-3}{k-2}
=
\binom{p+k-3}{p-2}+\binom{p+k-2}{p-2}
;
$$
see, e.g., Corollary~1.1.4 in \cite{Dai2013}.
Spherical harmonics associated with different values of~$k$ are orthogonal: if~$\varphi\in \mathcal{H}^p_k$ and~$\psi\in \mathcal{H}^p_\ell$ with~$k\neq \ell$, then~$\langle \varphi,\psi\rangle=0$. Therefore, merging orthonormal bases of several harmonic spaces~$\mathcal{H}^p_k$ provides an orthonormal basis of the direct sum of these spaces. Moreover, the collection of all spherical harmonics is dense in~$L^2(\mathcal{S}^{p-1},\nu_p)$, which coincides with the direct sum of~$\mathcal{H}^p_k$, $k=0,1,\ldots$; see, e.g., Theorem~2.2.2 in \cite{Dai2013}.

For each value of~$k$, we now provide a specific orthonormal basis of~$\mathcal{H}^p_k$ that is convenient for computations. For~$k=0$, we have~$d_{p,0}=1$ and a corresponding basis of~$\mathcal{H}^p_k$ is~$\{g_{1,0}\}$, with~$g_{1,0}(\xb)=1$ for any~$\xb$. We may therefore restrict to positive values of~$k$. For $p=2$, we have~$d_{2,k}=2$ for any positive integer~$k$, and the functions that are defined using polar coordinates~$\xb=(x_1,x_2)'=(\cos\theta, \sin\theta)' \in \mathcal{S}^{1}$ (with~$0\leq \theta<2\pi$) through
\begin{align}
	\label{basisp2}
	g_{1,k}(\xb)=\sqrt{2} \cos(k \theta)
	\quad
	\textrm{ and }
	\quad
	g_{2,k}(\xb)=\sqrt{2} \sin(k \theta)
\end{align}
form an orthonormal basis of~$\mathcal{H}^p_k$. The case~$p\geq 3$ is substantially more intricate. Consider the hyperspherical coordinates where~$\xb=(x_1,\ldots,x_p)' \in \mathcal{S}^{p-1}$ writes
\begin{align*}
	\left\{
	\begin{aligned}
		x_1 =&\; \sin \theta_{p-1} \cdots \sin \theta_2 \sin \theta_1,\\
		x_2 =&\; \sin \theta_{p-1} \cdots \sin \theta_2 \cos \theta_1,\\
		&\vdots\\
		x_{p-1} =&\; \sin \theta_{p-1} \cos \theta_{p-2},\\
		x_{p} =&\; \cos \theta_{p-1},
	\end{aligned}
	\right.
\end{align*}
with~$0\leq \theta_1<2\pi$ and~$0\leq \theta_j\leq \pi$ for~$j=2,\ldots,p-1$. Let further~$C_q^{\lambda}$ be the Gegenbauer polynomial of degree~$q$, defined for any~$\lambda> 0$ and~$q\in\N_0$, through
\begin{align}
	\label{Gegendefin}
	C_q^{\lambda}(t)
	:=
	\frac{(\lambda)_q 2^q}{q!}
	t^q
	\,
	{}_2 F_1({\textstyle{-\frac{q}{2},\frac{1-q}{2};1-q-\lambda;\frac{1}{t^2}}})
	=
	\sum_{j=0}^{\lfloor q/2 \rfloor}  (-1)^{j} c^\lambda_{q,j} t^{q-2j}
	,
\end{align}
where~$_2 F_1$ is the Gaussian hypergeometric function and where
$$
c^\lambda_{q,j}
:=
\frac{2^{q-2j}\Gamma(q-j+\lambda)}{\Gamma(\lambda) j! (q-2j)!}.
$$
Gegenbauer polynomials form an orthogonal basis on $L^2_\lambda([-1,1])$, the space of square-integrable functions with respect to the weight $x\mapsto (1-x^2)^{\lambda-1/2}$.
For~$p=3$, Gegenbauer polynomials reduce to Legendre polynomials ($\lambda=1/2$ in \eqref{Gegendefin}). In the sequel, we will also need to consider the Chebyshev polynomial of the first kind with degree~$q$, which is defined through~$C_q^{0}(\cos \theta):=\cos(q \theta)$, hence is given by
\begin{align}
	\label{Chebydefin}
	C_q^{0}(t)
	=
	\sum_{j=0}^{\lfloor q/2 \rfloor}  (-1)^{j} c^{0}_{q,j} t^{q-2j}
	,
	\ \
	\textrm{with }
	\
	c^{0}_{q,j}
	:=
	\frac{2^{q-2j-1} q \Gamma(q-j)}{j! (q-2j)!}
	\cdot
\end{align}
Chebyshev and Gegenbauer polynomials are connected by $\lim_{\lambda\to0}\lambda^{-1}C_k^{\lambda}(s)=(2/k)C_k^{0}(s)$, $s\in[-1,1]$. The following result provides an orthonormal basis of~$\mathcal{H}^p_k$ for~$p\geq 3$ and~$k>0$; see Theorem~1.5.1 in \cite{Dai2013}. The proposition incorporates corrections to several errata present in the statement of such Theorem~1.5.1.

\begin{proposition}
	\label{gigo}
	Fix an integer~$p \geq 3$ and a positive integer~$k$. Let~$\mathcal{M}_k=\{{\bf m}
	\in\N_0^p: |{\bf m}|:=m_1+\cdots+m_p=k \textrm{ and }m_p\in\{0,1\}\}$. For any~${\bf m}\in\mathcal{M}_k$, let
	$$
	\zeta_{\bf m}(s)
	:=
	\begin{cases}
		\cos(m_{p-1}s) & \text{if } m_p=0,  \\
		\sin((m_{p-1}+1)s) & \textrm{if } m_p=1
	\end{cases}
	\text{ and }
	B_{\bf m}
	:=
	b_{\bf m}
	\prod_{j=1}^{p-2}
	\frac{m_j!(\frac{p-j+1}{2})_{|{\bf m}^{j+1}|}(m_j+\lambda_j)}{(2\lambda_j)_{m_j}(\frac{p-j}{2})_{|{\bf m}^{j+1}|}\lambda_j}
	,
	$$
	where we write~$|{\bf m}^{j}|=m_j+\cdots+m_{p}$ and~$\lambda_j:=|{\bf m}^{j+1}|+(p-j-1)/2$ for any~$j=1,\ldots,p-2$, and where~$b_{\bf m}:=2$ if~$m_{p-1}+m_p>0$ and~$b_{\bf m}:=1$ otherwise.
	Then, defining~$\varphi_{\bf m}:\mathcal{S}^{p-1}\to\R$ in spherical coordinates through
	\begin{align}
		\label{grks}
		\varphi_{\bf m}(\xb)
		:=
		\sqrt{B_{\bf m}}
		\,
		r^{|{\bf m}|}
		\zeta_{\bf m}(\theta_1)
		\prod_{j=1}^{p-2}
		(\sin \theta_{p-j})^{|{\bf m}^{j+1}|}
		C_{m_j}^{\lambda_j}(\cos \theta_{p-j})
		,
	\end{align}
	the collection~$\{\varphi_{\bf m}:{\bf m}\in\mathcal{M}_k\}$ is a real orthonormal basis of~$\mathcal{H}^p_k$, in the sense that~$\langle \varphi_{\bf m},\varphi_{\tilde{\bf m}}\rangle=\delta_{{\bf m},\tilde{{\bf m}}}$ for any~${\bf m},\tilde{{\bf m}}\in\mathcal{M}_k$.
\end{proposition}

Note that this orthonormal basis is indeed of cardinality~$d_{p,k}$. To work with a specific basis in the rest of the paper, we adopt the following notation. Let us enumerate the~${\bf m}$'s in~$\mathcal{M}_k$ as~${\bf m}_r$, $r=1,\ldots,d_{p,k}$, in an arbitrary way but for the fact that~${\bf m}_1=(k,0,\ldots,0)\pr$ (the exact ordering for~$r=2,\ldots,d_{p,k}$ will not be important in the sequel). We then let
\begin{align}
	\label{eq:Yk0}
	g_{1,k}(\xb)
	&=
	\varphi_{{\bf m}_1}(\xb)
	=
	\frac{1}{\sqrt{d_{p,k}}}
	\Big(1+ \frac{2k}{p-2}\Big)
	C_k^{(p-2)/2}(\cos \theta_{p-1})
	,
	\\%[2mm]
	g_{r,k}(\xb)
	&=
	\varphi_{{\bf m}_r}(\xb)
	,
	\quad
	r=2,\ldots,d_{p,k}
	.
	\label{eq:Ykr}
\end{align}
As shown in the later sections, $g_{1,k}$ is the only basis function of~$\mathcal{H}^p_k$ that will play an essential role in the paper. The exact expression of the remaining basis functions in~\eqref{eq:Ykr}, which is rather cumbersome (see~\eqref{grks}), will play a very minor role only.

Irrespective of~$p$, one has the \emph{addition formula}
\begin{align}
	\label{iin}
	\sum_{r=1}^{d_{p,k}}
	g_{r,k}(\ub)g_{r,k}(\vb)
	=
	\begin{cases}
		2 \cos(k \angle(\ub,\vb)) = 2 C_k^{0}(\ub \pr \vb) & {\rm for} \; p=2, \\
		\big(1+ \frac{2k}{p-2}\big)
		C_k^{(p-2)/2} ({\bf u}\pr \vb) & {\rm for} \; p\geq 3
		,
	\end{cases}
\end{align}
where~$\angle(\ub,\vb):= \arccos (\ub\pr \vb)$ is the angle between~$\ub,\vb\in\mathcal{S}^{p-1}$; see Equation~(1.2.8) in \cite{Dai2013} for~$p\geq 3$. The name follows from the fact that this result extends to~$p\geq 3$ what is the usual addition formula of the cosine function for~$p=2$. %

Now, for any
real sequence~$v=(v_k)$ such that~$\sum_{k=1}^\infty v_k^2 d_{p,k}<\infty$, consider the mapping~${\rm t}:\mathcal{S}^{p-1} \to L^2(\mathcal{S}^{p-1},\nu_p)$ that maps~$\ub$ to~${\rm t}(\ub):=\sum_{k=1}^\infty v_k {\rm t}_k(\ub)$,
where we let
${\rm t}_k(\ub):=\sum_{r=1}^{d_{p,k}} g_{r,k}(\ub)g_{r,k}$ for any~$ \ub \in\mathcal{S}^{p-1}$; see~\eqref{basisp2} for~$p=2$ and the discussion below Proposition~\ref{gigo} for~$p\geq 3$. Considering the mapping ${\bf G}_{p,k}:\mathcal{S}^{p-1}\to\mathbb{R}^{d_{p,k}}$ defined as~${\bf G}_{p,k}:=(g_{1,k},\ldots,g_{d_{p,k},k})'$ %
allows us to formally define the Sobolev statistic (in the sequel we put $\Ub_i:=\Ub_i (\Gamb)$)
\begin{align}
	S^{(n)}_v
	:=&\; \frac{1}{n} \sum_{i,j=1}^n \langle {\rm t}(\Ub_i),  {\rm t}(\Ub_j) \rangle =\frac{1}{n} \sum_{i,j=1}^n\sum_{k=1}^\infty v_k^2\langle {\rm t}_k(\Ub_i),  {\rm t}_k(\Ub_j) \rangle
	\nonumber\\
	=&\;\frac{1}{n} \sum_{i,j=1}^n\sum_{k=1}^\infty v_k^2
	\bigg(
	\sum_{r=1}^{d_{p,k}}
	g_{r,k}(\Ub_i)g_{r,k}(\Ub_j)
	\bigg)
	= \sum_{k=1}^\infty v_k^2 \bigg\| \frac{1}{\sqrt{n}} \sum_{i=1}^n {\bf G}_{p,k}(\Ub_i) \bigg\|^2 ,\!\!\!\!
	\label{Sobolgen}
\end{align}
where we used the fact that, for any~$\ub,\vb\in\mathcal{S}^{p-1}$, we have~$\langle {\rm t}_k(\ub),{\rm t}_\ell(\vb)\rangle=0$ for~$k\neq \ell$ and~$\langle{\rm t}_k(\ub), {\rm t}_k(\vb) \rangle= \sum_{r=1}^{d_{p,k}} g_{r,k}(\ub)g_{r,k}(\vb)=({\bf G}_{p,k}(\ub))\pr {\bf G}_{p,k}(\vb)$. The addition formula above now allows us to write~$S^{(n)}_v$ as in~\eqref{Sobol1}, which is the formula that is used in practice to evaluate this test statistic. Since the expression of~$S^{(n)}_v$ in~\eqref{Sobol1} involves the observations~$\Ub_i$, $i=1,\ldots,n$, only through the inner products~$\Ub_i'\Ub_j$, Sobolev tests are rotation-invariant in the sense that~$S^{(n)}_v(\Ob\Ub_1,\ldots,\Ob\Ub_n)=S^{(n)}_v(\Ub_1,\ldots,\Ub_n)$ for any $p\times p$ orthogonal matrix~$\Ob$.

We now describe the asymptotic null behavior of Sobolev test statistics. Write~${\rm P}_0\n$ for the hypothesis under which the observations~$\Ub_1,\ldots,\Ub_n$ form a random sample from the uniform distribution on~$\mathcal{S}^{p-1}$. Throughout, expectations and variances under~${\rm P}_0\n$ will be denoted as~${\rm E}_{0}\n[\cdot]$ and~${\rm Var}_{0}\n[\cdot]$, respectively. Assume for a moment that the Sobolev test statistic~$S^{(n)}_v$ is finite, with~$K_v:=\max\{k:v_k\neq 0\}$ say. Then, the last expression in~\eqref{Sobolgen} entails that
\begin{align}
	\label{truncSobol}
	S^{(n)}_v
	=
	n
	({\bf T}_{p,K_v}\n)\pr {\rm diag}(v_1^2 {\bf I}_{d_{p,1}}, \ldots, v_{K_v}^2 {\bf I}_{d_{p,K_v}} )
	{\bf T}_{p,K_v}\n
	,
\end{align}
with
$$
{\bf T}_{p,K}\n
:=
\frac{1}{n}
\sum_{i=1}^n
\left( ({\bf G}_{p,1}(\Ub_i))', \ldots, ({\bf G}_{p,K}(\Ub_i))' \right)\pr
$$
a vector of length $d_{p,1}+\cdots+d_{p,K}$. For any positive integer~$k$, the orthonormality of the functions~$\{g_{r,k}\}_{r=1}^{d_{p,k}}$ and mutual orthogonality of the subspaces~$\mathcal{H}^p_k$, $k=0,1,\ldots$, entail that
\begin{align}
	\label{MomGnull1}
	{\rm E}_0\n[{\bf G}_{p,k}(\Ub_1)]&={\bf 0}\quad\text{and} \quad {\rm Var}_0\n[{\bf G}_{p,k}(\Ub_1)]={\bf I}_{d_{p,k}},\\
	\label{MomGnull3}
	{\rm E}_0\n[{\bf G}_{p,k}(\Ub_1)({\bf G}_{p,\ell}(\Ub_1))']&={\bf 0},
	\quad
	k\neq \ell
	.
\end{align}
The central limit theorem readily implies that $S^{(n)}_v$ converges weakly under~${\rm P}_0\n$ to~$\sum_{k=1}^{K_v} v_k^2 \chi^2_{d_{p,k}}$. %
The result actually extends to infinite Sobolev tests, under the summability condition provided in the following theorem; see Theorem~4.1 and Equation~(4.2) in \cite{Gine1975}.

\begin{theorem}
	\label{asymptnull}
	Fix an integer~$p\geq 2$ and a sequence~$v=(v_k)$ such that~$\sum_{k=1}^{\infty}
	v_k^2
	d_{p,k}
	< \infty
	$.
	Then, under~${\rm P}\n_{0}$,
	$
	S^{(n)}_v
	\stackrel{\mathcal{D}}{\to}
	\sum_{k=1}^{\infty} v_k^2 \chi^2_{d_{p,k}}
	$.
\end{theorem}

The resulting infinite Sobolev tests reject the null hypothesis of uniformity at asymptotic level~$\alpha$ when~$S^{(n)}_v$ exceeds the upper $\alpha$-quantile of the distribution of~$\sum_{k=1}^{\infty} v_k^2 \chi^2_{d_{p,k}}$.

%-------------------------------%
\section{Asymptotics under the alternative: non-null moments of random spherical harmonics}
\label{sec:moments}
%-------------------------------%

In the sequel, our goal is to provide an extensive study of the local asymptotic behavior of Sobolev tests of uniformity of the $\Ub_i=\Ub_i(\Gamb)$'s under the most commonly considered alternatives, namely under (absolutely continuous) rotationally symmetric ones. In other words, we consider alternatives associated with densities of the form
\begin{align}
	\label{Rotsymfirst}
	\ub
	\mapsto c_{p,\kappa, f} f(\kappa \ub\pr \thetab)
	,
\end{align}
where $\thetab \in \mathcal{S}^{p-1}$ is a location parameter, $\kappa>0$ is a concentration parameter, $f:\R\to\R^+$ is an \emph{angular function} with~$f(0)=1$ (for standardization), and
\begin{align}
	\label{defcpkf}
	c_{p,\kappa, f}:= 1 \, \Big/\int_{-1}^1 (1-s^2)^{(p-3)/2} f(\kappa s) \,\mathrm{d}s
\end{align}
is a normalizing constant. These parameters are not identified~(the concentration parameter~$\kappa$ may always be included in~$f$; also, when~$f$ is symmetric about zero, $\thetab$ is only identified up to a sign), but this will not be an issue in the sequel since these distributions are only providing alternatives under which to investigate the asymptotic properties of perfectly well-defined tests. Many well-known distributions on the sphere fit in the setup of~\eqref{Rotsymfirst}. For instance, the von Mises--Fisher distributions introduced above, which are often regarded as the $\mathcal{S}^{p-1}$-equivalent of isotropic Gaussian distributions on~${\mathbb R}^p$, are obtained with~$f(s)=\exp(s)$. Watson distributions \citep{Watson1965}, that are associated with~$f(s)=\exp(s^2)$, provide a single-spiked version of the \cite{Bingham1974} distributions. Further distributions of this type are the directional Cauchy distributions with angular function~$f(s)=(1+2s)^{-1}$ \citep{Garcia-Portugues2013a} and the exponential distributions with angular functions~$f(s)= \exp(s^b)$ for $b>0$ \citep{Paindaveine2020}, as well as other distributions that can be reparametrized to fit \eqref{Rotsymfirst} (see, e.g., \citealt{McCullagh1989}, \citealt{Jones2005}, or \citealt{Kato2020}). 

We investigate here the limiting behavior of Sobolev tests under local alternatives of the form~\eqref{Rotsymfirst}, that is, when~$\kappa$ is replaced in~\eqref{Rotsymfirst} with a sequence~$(\kappa_n)$ converging to zero at a suitable rate. Accordingly, write~${\rm P}_{\thetab,\kappa_n,f}\n$ for the hypothesis under which the observations~$\Ub_1,\ldots,\Ub_n$ form a random sample from a rotationally symmetric distribution with location~$\thetab$, concentration~$\kappa_n$, and angular function~$f$; see~\eqref{Rotsymfirst}. From rotation-invariance, the power of Sobolev tests under~${\rm P}_{\thetab,\kappa_n,f}\n$ will not depend on~$\thetab$, so we may safely restrict in the sequel to the case~$\thetab=\thetab_0=(0,\ldots,0,1)'\in\R^p$. The resulting sequence of hypotheses will be denoted as~${\rm P}_{\kappa_n,f}\n$, and the corresponding expectation and variance will be denoted as~${\rm E}_{\kappa_n, f}\n[\cdot]$ and~${\rm Var}\n_{\kappa_n, f}[\cdot]$, respectively. Note that
\begin{enumerate}[label=(\textit{\roman{*}}),ref=\textit{\roman{*}}]
	\item for $\Gamb= {\bf I}_d$ (the classical uniformity problem), the alternatives above are classical rotationally symmetric alternatives to uniformity on $\mathcal{S}^{d-1}$, while
	\item taking $\Gamb$ such that $\Gamb\pr \Gamb= {\bf I}_d- {\pmb \mu}_1 {\pmb \mu_1}\pr$ for some ${\pmb \mu}_1 \in \mathcal{S}^{d-1}$ (the problem of testing for rotational symmetry), the alternatives above extend the \emph{tangent von Mises alternatives} and \emph{tangent elliptical alternatives} introduced in \cite{Garcia-Portugues2020}.
\end{enumerate}
The discussion at the end of Section \ref{sec:Sobol} makes it clear that the asymptotic non-null behavior of Sobolev tests will be governed by the non-null counterparts of~\eqref{MomGnull1}--\eqref{MomGnull3}. %

\begin{proposition}
	\label{Exprop}
	Fix an integer~$p\geq 2$ and a positive integer~$k$. Let~$(\kappa_n)$ be a positive real sequence that is~$o(1)$ and~$f$ be an angular function that is continuous at zero. Then, as~$n$ diverges to infinity,
	\begin{align}
		\label{shiftstep1}
		{\rm E}\n_{\kappa_n,f}[{\bf G}_{p,k}(\Ub_1)]
		&=
		t_{p,k}
		\,
		\mathrm{E}_{\kappa_n,f}\n\big[C_k^{(p-2)/2}(\Ub_1'\thetab_0)\big]
		{\bf e}_{1}
		,\\
		\label{MomGnonnull1}
		{\rm Var}_{\kappa_n,f}\n[{\bf G}_{p,k}(\Ub_1)]&={\bf I}_{d_{p,k}}
		+o(1)
		,\\
		\label{MomGnonnull2}
		{\rm E}_{\kappa_n,f}\n[{\bf G}_{p,k}(\Ub_1)({\bf G}_{p,\ell}(\Ub_1))']&=o(1)
		,
		\quad
		k\neq \ell,
	\end{align}
	with~$t_{p,k}:=\sqrt{2}$ for~$p=2$ and $t_{p,k}:=
	(1+2k/(p-2))/\sqrt{d_{p,k}}$ for~$p\geq 3$.
\end{proposition}

Comparing~\eqref{shiftstep1}--\eqref{MomGnonnull2} with~\eqref{MomGnull1}--\eqref{MomGnull3}, this result shows that the behavior of Sobolev tests under local alternatives will be exclusively driven by the shift the expectations in~\eqref{shiftstep1} show compared to the null hypothesis. In order to derive the corresponding limiting powers, we need to describe the way the expectation in the right hand side of \eqref{shiftstep1} behaves as a function of~$\kappa_n$. This expectation involves Chebyshev %
or Gengenbauer polynomials, %
hence\nopagebreak[4] in all cases can be decomposed in terms of the non-null moments~$e_{n,m}:= {\rm E}\n_{\kappa_n,f}[(\Ub_1'\thetab_0)^m]$, $m=0,1,\ldots$ We recall that the corresponding null moments are given by
\begin{align}
	a_m
	:=&\;
	{\rm E}\n_0[(\Ub_1'\thetab_0)^m]
	=
	c_{p}
	\int_{-1}^1 s^m (1-s^2)^{(p-3)/2}  \,\mathrm{d}s
	\nonumber\\
	=&\;
	\mathbb{I}[{\textstyle{\frac{m}{2}}} \in {\mathbb N} ]
	\frac{\Gamma(\frac{m+1}{2})\Gamma(\frac{p}{2})}{\sqrt{\pi}\Gamma(\frac{m+p}{2})}
	=
	\mathbb{I}[{\textstyle{\frac{m}{2}}} \in {\mathbb N} ]
	\prod_{r=0}^{{\textstyle{\frac{m}{2}}}-1}
	\frac{1+2r}{p+2r}
	,
	\quad
	m=0,1,\ldots,
	\label{foref1}
\end{align}
with the convention that a product involving no factor is equal to one; see, e.g., Lemma~A.1 in \cite{Paindaveine2016}. The next result then describes how the non-null moments~$e_{n,m}$ behave in terms of~$\kappa_n$.

\begin{proposition}
	\label{Momentprop}
	Fix an integer~$p\geq 2$ and integers $q\geq m>0$. Let~$(\kappa_n)$ be a positive real sequence that is~$o(1)$ and~$f$ be an angular function that is~$q-m$ times differentiable at zero. Define the coefficients~$b_{m,i}$ through
	$
	(b_{m,0},b_{m,1},\ldots,b_{m,q-m})'
	:=
	\Ab_{q-m}^{-1}\vb_{q-m}^{(m)}
	,
	$
	where the $(t+1)\times (t+1)$ matrix~$\Ab_{t}$ and $(t+1)$-vector~$\vb_t^{(m)}$ are defined as
	$
	\Ab_{t}
	=
	\big(
	a_{i-j} f^{\underline{i-j}}(0)/(i-j)!\;  \mathbb{I}[i \geq j]
	\big)_{i,j=0,1,\ldots,t}$
	and
	$\vb_t^{(m)}
	=
	\big(
	a_{m+i} f^{\underline{i}}(0)/i!
	\big)_{i=0,1,\ldots,t}
	.
	$
	Then, as~$n$ diverges to infinity,
	$$
	e_{n,m}
	=
	\sum_{\ell=0}^{q-m}
	b_{m,\ell}
	\kappa_n^{\ell}
	+o({\kappa_n^{q-m}})
	$$
\end{proposition}

This result, that is of independent interest, is an important step to derive the asymptotic behavior of the expectation in the right-hand side of~\eqref{shiftstep1}. To state the result, we need the following further notation. For any positive integer~$i$, define the coefficients~$m_{k,i}$, $k=0,1,\ldots,i$, through
\begin{align} \label{monomials}
	\sum_{k=0}^i m_{k,i} C_{k}^{(p-2)/2}(t)
	=
	t^{i}.
\end{align}
These coefficients are uniquely defined since %
$\{C_{k}^{(p-2)/2}\}_{k=0}^i$, form a basis of the vector space of polynomials of degree at most~$i$. Note that, irrespective of~$i$, we have~$m_{i-\ell,i}=0$ for $\ell$~odd (see Equation~(2.29) in \citealp{Wunsche2017}), so that a necessary condition for~$m_{k,i}\neq 0$ is that~$k\sim i$ (i.e.,  that~$k$ and~$i$ have the same parity). We then have the following result.

\begin{proposition}
	\label{finalExprop}
	Fix integers~$p\geq 2$, $k>0$, and~$r\geq 0$. Let~$(\kappa_n)$ be a positive real sequence that is~$o(1)$ and~$f$ be an angular function that is $k+r$ times differentiable at zero. Let ${\bf z}^{(k)}_{t}$ be the $(t+1)$-vector
	$
	{\bf z}_t^{(k)}
	=
	\big(
	m_{k,i} f^{\underline{i}}(0)/i!\; \mathbb{I}[k\leq i]
	\big)_{i=0,1,\ldots,t}
	$
	and~${\bf a}_{i;k+r}\pr$ be the $i$th row of~${\bf A}^{-1}_{k+r}$.
	Then, as~$n$ diverges to infinity,
	$$
	\mathrm{E}\n_{\kappa_n,f} [C_k^{(p-2)/2}(\Ub_1'\thetab_0)]
	=
	\frac{1}{t_{p,k}^2}
	\sum_{\ell=k}^{k+r}
	\,
	({\bf a}_{\ell+1;k+r}\pr {\bf z}^{(k)}_{k+r})
	\kappa_n^\ell+o(\kappa_n^{k+r})
	$$
\end{proposition}

Note that~$\Ab_{t}=\mathbf{I}_{t+1}+\Lb_t$, where~$\Lb_t$ is a strictly lower triangular matrix. Since~$\Lb_t^{t+1}={\bf 0}$, we have
${\bf A}_t^{-1}= \sum_{j=0}^{t} (-1)^j \Lb_t^j$,
which provides an explicit expression for the $(k+r+1)$-vectors~${\bf a}_{\ell+1;k+r}$ involved in Proposition~\ref{finalExprop}. In particular, this shows that this vector is always of the form~${\bf a}_{\ell+1;k+r}=(\lambda_1,\ldots,\lambda_\ell,1,0,\ldots,0)'$. Since~$a_m=0$ for any odd positive integer~$m$, the entries~$\lambda_{j}$ with~$j\sim \ell$ are actually equal to zero.
Now, taking~$r=0$ in Proposition~\ref{finalExprop} yields
$$
\mathrm{E}\n_{\kappa_n,f}\big[C_k^{(p-2)/2}(\Ub_1'\thetab_0)\big]
=
\frac{1}{t_{p,k}^2}
({\bf a}_{k+1;k}\pr {\bf z}^{(k)}_{k})
\kappa_n^k+o(\kappa_n^{k})
.
$$
From the discussion on the structure of~${\bf a}_{\ell+1;k+r}$ above and the fact that the first~$k$ components of the $(k+1)$-vector~${\bf z}_k^{(k)}
$ are equal to zero, this rewrites
$$
\mathrm{E}\n_{\kappa_n,f}\big[C_k^{(p-2)/2}(\Ub_1'\thetab_0)\big]
=
\frac{{\bf e}_{k+1}' {\bf z}_k^{(k)}}{t_{p,k}^2}
\kappa_n^k+o(\kappa_n^{k})
=
\frac{m_{k,k} f^{\underline{k}}(0)}{(k!)t_{p,k}^2}
\kappa_n^k+o(\kappa_n^{k})
.
$$
Using Proposition~\ref{Exprop}, this entails that if the $k$th derivative of~$f$ at zero exists, then
$$
\mathrm{E}_{\kappa_n,f}\n[{\bf G}_{p,k}(\Ub_1)]
=
\frac{m_{k,k} f^{\underline{k}}(0)}{(k!)t_{p,k}}
\kappa_n^k
{\bf e}_{1}
+
o(\kappa_n^{k})
.
$$
This type of results will be most important when computing the local powers of Sobolev tests in the next sections.

%-------------------------------%
\section{Asymptotics under the alternative: local non-null behavior of Sobolev tests}
\label{sec:localpowfiniteinf}
%-------------------------------%

We now determine the asymptotic behavior of finite Sobolev tests under the local alternatives considered in the previous sections. We first do so for finite tests, placing minimal differentiability assumptions on the angular function~$f$ at zero (Theorem \ref{Sobolfin}), and then extend these results to infinite tests assuming stronger differentiability at zero (Theorem~\ref{Sobolfininfin}).

\subsection{Finite Sobolev tests}

Recall that a Sobolev test is finite if finitely many of the coefficients~$v_k$ defining this test are non-zero. Letting~$K_v:=\max\{k:v_k\neq 0\}$ throughout this section, finite Sobolev tests are thus based on statistics of the form
\begin{align}
	\label{truncSobolagain}
	S^{(n)}_v
	=
	n
	({\bf T}_{p,K_v}\n)\pr {\rm diag}(v_1^2 {\bf I}_{d_{p,1}}, \ldots, v_{K_v}^2 {\bf I}_{d_{p,K_v}} )
	{\bf T}_{p,K_v}\n
	.
\end{align}
As we saw in Section~\ref{sec:Sobol}, the asymptotic (standard) normality of~$\sqrt{n}{\bf T}_{p,K_v}\n$ under the null hypothesis of uniformity readily follows from \eqref{MomGnull1}--\eqref{MomGnull3} and from the central limit theorem, which implies that, still under the null hypothesis,
$$
S^{(n)}_v
\stackrel{\mathcal{D}}{\to}
\sum_{k=1}^{K_v} v_k^2 \chi^2_{d_{p,k}}.
$$
We have the following result under the local alternatives from the previous sections.

\begin{proposition}
	\label{asymptnormtrunc}
	Fix integers~$p\geq 2$ and~$K>0$. Let~$(\kappa_n)$ be a positive real sequence that is~$o(1)$ and~$f$ be an angular function that is continuous at zero. Then,
	$
	\sqrt{n}
	\big(
	{\bf T}_{p,K}\n- {\rm E}\n_{\kappa_n, f}[{\bf T}_{p,K}\n]
	\big)
	$
	is asymptotically standard normal under~${\rm P}_{\kappa_n,f}\n$.
\end{proposition}

The asymptotic power of Sobolev tests based on statistics of the form~\eqref{truncSobolagain} will then follow by combining Proposition~\ref{asymptnormtrunc} with Propositions~\ref{Exprop} and~\ref{finalExprop}. %

\begin{theorem}
	\label{Sobolfin}
	Fix an integer~$p\geq 2$ and a sequence~$v=(v_k)$ with only finitely many non-zero terms. Let~$f$ be an angular function that is~$K_v$ times differentiable at zero. Let~$k_v:=\min\{k:v_k\neq 0\}\leq K_v$, fix~$\tau>0$, and put
	\begin{align}
		\label{ksith41}
		\xi_{p,k}(\tau)
		:=
		\frac{m^2_{k,k}(f^{\underline{k}}(0))^2 \tau^{2k}}{(k!)^2 t_{p,k}^2}
		=
		\frac{d_{p,k} (f^{\underline{k}}(0))^2
			\tau^{2k}}{\prod_{\ell=0}^{k-1} (p+2\ell)^2}
		\cdot
	\end{align}
	Then, %
	\begin{enumerate}[label=(\textit{\roman{*}}),ref=\textit{\roman{*}}]
		\item under~${\rm P}\n_{\kappa_n, f}$, with $\kappa_n=o(n^{-1/(2{k_v})})$,
		$
		S^{(n)}_v
		\stackrel{\mathcal{D}}{\to}
		\sum_{k=1}^{K_v} v_k^2 \chi^2_{d_{p,k}}
		$; \label{Sobolfin:1}
		\item under~${\rm P}\n_{\kappa_n, f}$, with $\kappa_n=n^{-1/(2{k_v})} \tau$,
		$
		S^{(n)}_v
		\stackrel{\mathcal{D}}{\to}
		\sum_{k=1}^{K_v} v_k^2 \chi^2_{d_{p,k}}(\delta_{k,k_v}\xi_{p,k_v}(\tau))
		$;\label{Sobolfin:2}
		\item under~${\rm P}\n_{\kappa_n, f}$, where~$n^{1/(2{k_v})}\kappa_n\to\infty$ and~$f^{\underline{k_v}}(0)\neq 0$,
		$
		{\rm P}\n_{\kappa_n,f}[S^{(n)}_v>c]\to 1
		$,
		so that the corresponding Sobolev test is consistent under this sequence of alternatives. \label{Sobolfin:3}
	\end{enumerate}
\end{theorem}

This result shows that the consistency rate of a finite Sobolev test is typically~$\kappa_n\sim n^{-1/(2{k_v})}$. More precisely, provided that~$f^{\underline{k_v}}(0)\neq 0$, such a test will show non-trivial asymptotic powers (that is, asymptotic powers strictly between the nominal level~$\alpha$ and one) against local alternatives of the form~$\kappa_n=n^{-1/(2{k_v})}\tau$, with~$\tau>0$, and will be blind (resp., consistent) against less severe (resp., more severe) alternatives. Using the terminology from \cite{Bha2019a,Bha2019}, $\kappa_n\sim n^{-1/(2{k_v})}$ is thus the \emph{detection threshold} of a finite Sobolev test when~$f^{\underline{k_v}}(0)\neq 0$. In the case where~$f^{\underline{k_v}}(0)=0$, a finite Sobolev test will be blind to local alternatives in~$\kappa_n=O(n^{-1/(2{k_v})})$, but Theorem~\ref{Sobolfin} remains silent about the possible existence of more severe alternatives that could be detected by such a test. We will address this point in Theorem~\ref{Sobolfinbisbis} below.

While Theorem~\ref{Sobolfin} applies to any finite Sobolev test, let us consider two specific examples. We first focus on the Rayleigh test (obtained with~$v_k=\delta_{k,1}$). Theorem~\ref{Sobolfin} entails that in the von Mises--Fisher case $f(s)=\exp(s)$ or, more generally, under any angular function with~$f'(0)\neq 0$, this test shows non-trivial asymptotic powers against the (contiguous) alternatives associated with~$\kappa_n\sim n^{-1/2}$. More precisely, Theorem~\ref{Sobolfin} entails that, under~${\rm P}\n_{\kappa_n, f}$ with $\kappa_n=n^{-1/2} \tau$, the Rayleigh test statistic is asymptotically
$$
\chi^2_{d_{p,1}} \bigg( \frac{1}{p^2}d_{p,1} (f'(0))^2 \tau^2 \bigg)
=
\chi^2_{p} \bigg( \frac{1}{p} (f'(0))^2 \tau^2 \bigg)
,
$$
which agrees with Equation~(3.4) in \cite{Cutting2017} (note that~$f'(0)=1$ there). Theorem~\ref{Sobolfin} also allows us to conclude that, in the Bingham case $f(s)=\exp(s^2)$, for which~$f'(0)=0$, the Rayleigh test is blind to the local alternatives associated with~$\kappa_n\sim n^{-1/2}$. Whether or not it can detect more severe alternatives is thus the issue we will consider below in Theorem~\ref{Sobolfinbisbis}.
We then turn to the Bingham test (obtained with~$v_k=\delta_{k,2}$). Theorem~\ref{Sobolfin} entails that, both in the von Mises--Fisher and Bingham cases, this test shows non-trivial asymptotic powers against local alternatives associated with~$\kappa_n\sim n^{-1/4}$.
More precisely, in the Bingham case~$f(s)=\exp(s^2)$ or, more generally, under any angular function with~$f''(0)\neq 0$ (including the von Mises--Fisher ones), the Bingham test statistic, under~${\rm P}\n_{\kappa_n, f}$ with~$\kappa_n=n^{-1/4} \tau$, is asymptotically
$$
\chi^2_{d_{p,2}} \bigg(\frac{d_{p,2} (f''(0))^2 \tau^{4}}{p^2(p+2)^2}\bigg)
=
\chi^2_{(p-1)(p+2)/2} \bigg(\frac{(p-1)(f^{\prime \prime}(0))^2 \tau^4}{2 p^2(p+2)} \bigg)
,
$$
which is compatible with Equation~(9) in~\cite{Cutting2020} (there, the Bingham test shows non-trivial asymptotic powers against Bingham alternatives in~$\kappa_n\sim n^{-1/2}$ rather than~$\kappa_n\sim n^{-1/4}$, but this is only due to the different choice of parametrization in~$\kappa$ adopted in that paper). The fact that the Bingham test detects von Mises--Fisher alternatives in~$\kappa_n\sim n^{-1/4}$ is also in line with Equation~(14) in \cite{Cutting2021}, although the latter result focuses on a high-dimensional asymptotic scenario under which~$p=p_n$ diverges to infinity as~$n$ does.

As explained above, Theorem~\ref{Sobolfin} fails to thoroughly describe the asymptotic non-null properties of finite Sobolev tests in the case where~$f^{\underline{k_v}}(0)=0$. In such a case, the result states that a finite Sobolev test is blind to alternatives in~$\kappa_n=O(n^{-1/(2k_v)})$ but it remains silent on what happens for more severe alternatives. This motivates the following result.
\begin{theorem}
	\label{Sobolfinbisbis}
	Fix an integer~$p\geq 2$ and a sequence~$v=(v_k)$ with only finitely many non-zero terms. Let~$f$ be an angular function that is~$q(\geq K_v)$ times differentiable at zero. Let~$k_v:=\min\{k:v_k\neq 0\}$ and fix~$\tau>0$. Then, %
	\begin{enumerate}[label=(\textit{\roman{*}}),ref=\textit{\roman{*}}]
		\item let~$k_*$ be the minimum value of~$k\in\{k_v,\ldots,q\}$ such that (a)~$f^{\underline{k}}(0)\neq 0$ and~(b) $\mathcal{V}_k:=\{\ell=k_v,\ldots,k: \ell
		\sim k
		\textrm{ and }v_\ell\neq 0\}$ is non-empty (assuming that such a~$k$ exists). Let~$k_\dagger:=\min\mathcal{V}_{k_*}$ and put
		\begin{align*}
			\xi_{p,k,k_*}(\tau)
			:=&\;
			\frac{m^2_{k,k_*}(f^{\underline{k_*}}(0))^2\tau^{2k_*}}{(k_*!)^2 t_{p,k}^2 }
		=%
		\frac{1}{(k_*!)^2}
		t_{p,k}^2  (f^{\underline{k}_*}(0))^2 \tau^{2k_*}
		\Bigg(
		\sum_{j=0}^{\lfloor k/2 \rfloor}
		(-1)^{j}
		c^{(p-2)/2}_{k,j}
		a_{k+k_*-2j}
		\Bigg)^2
		.
	\end{align*}
	Then, under~${\rm P}\n_{\kappa_n, f}$, with $\kappa_n=n^{-1/(2{k_*})} \tau$,
	$
	S^{(n)}_v
	\stackrel{\mathcal{D}}{\to}
	\sum_{k=1}^{K_v} v_k^2 \chi^2_{d_{p,k}}(\mathbb{I}[k_\dagger\leq k\leq k_*,\, k\sim k_*]\xi_{p,k,k_*}(\tau));\label{Sobolfinbisbis:1}
	$
	\item if no~$k\in\{1,\ldots,q\}$ do satisfy~(a)--(b), then, under ${\rm P}\n_{\kappa_n, f}$, with $\kappa_n= n^{-1/(2q)} \tau$,
	$
	S^{(n)}_v
	\stackrel{\mathcal{D}}{\to}
	\sum_{k=1}^{K_v} v_k^2 \chi^2_{d_{p,k}}
	$.\label{Sobolfinbisbis:2}%
\end{enumerate}
\end{theorem}

This result shows in particular that the ``pure'' Sobolev test associated with~$v_k=\delta_{k,\ell}$ will show non-trivial asymptotic powers against local alternatives in~$\kappa_n\sim n^{-1/(2{k_*})}$, where~$k_*$ is the smallest integer larger than or equal to~$\ell$, with the same parity as~$\ell$, and for which~$f^{\underline{k_*}}(0)\neq 0$ (provided, of course, that~$f$ is sufficiently many times differentiable at zero). The case~$k_*=\ell(=k_v)$ is covered by Theorem~\ref{Sobolfin}, but not the case~$k_*>\ell$. For instance, for alternatives associated with %
$f(s)=\exp(s^3)$ (for which~$f^{\underline{k}}(0)\neq 0$ if and only if~$k$ is a multiple of~3), the Rayleigh test shows non-trivial asymptotic powers against local alternatives in~$\kappa_n\sim n^{-1/6}$ ($k_*=3>1=\ell$). More precisely, Theorem~\ref{Sobolfinbisbis} entails that, under~${\rm P}\n_{\kappa_n, f}$ with $\kappa_n=n^{-1/6} \tau$, the Rayleigh test statistic is asymptotically
$$
\chi^2_{d_{p,1}} \bigg(   \frac{(f^{\underline{3}}(0))^2 }{4p(p+2)^2}
\tau^{6}
\bigg)
=
\chi^2_{p} \bigg(  \frac{9}{p(p+2)^2}
\tau^{6} \bigg)
$$
(in the notation of Theorem~\ref{Sobolfinbisbis}, we have $k_\dagger=1$ and $k_*=3$).

For the sake of simplicity, we did not include in Theorem~\ref{Sobolfinbisbis} results on detection thresholds in the sense described below Theorem~\ref{Sobolfin}. Inspection of the proof of Theorem~\ref{Sobolfinbisbis}\eqref{Sobolfinbisbis:1}, however, makes it clear that, provided that~$\xi_{p,k,k_*}(\tau)\neq 0$, the rate~$\kappa_n\sim n^{-1/(2{k_*})}$ is the detection threshold of the finite Sobolev test considered. Also, a direct consequence of Theorem~\ref{Sobolfinbisbis}\eqref{Sobolfinbisbis:2} is that Sobolev tests will be blind to local alternatives in~$\kappa_n\sim n^{-1/(2{k_*})}$ as soon as the non-zero coefficients~$(v_k)$ and the non-zero derivatives~$f^{\underline{k}}(0)$ show at ranks~$k$ with opposite parities. We have the following result.

\begin{corollary}
\label{Sobolfinbis}
Fix an integer~$p\geq 2$ and a sequence~$v=(v_k)$ with only finitely many non-zero terms and such that~$v_k=0$ for any~$k$ even (resp., odd). Let~$f$ be an angular function that is~$q(\geq K_v)$ times differentiable at zero and such that~$f^{\underline{k}}(0)=0$ for any~$k$ odd (resp., even) in~$\{k_v,\ldots,q\}$.
Then, under~${\rm P}\n_{\kappa_n, f}$, with $\kappa_n= n^{-1/(2q)} \tau$,
$
S^{(n)}_v
\stackrel{\mathcal{D}}{\to}
\sum_{k=1}^{K_v} v_k^2 \chi^2_{d_{p,k}}
.
$
\end{corollary}

In particular, applying this result with~$q$ arbitrarily large (assuming that~$f$ is infinitely many times differentiable at zero) shows that Sobolev tests with non-zero coefficients~$v_k$ at odd (resp., even) ranks only are blind---at any polynomial rate---to alternatives such that~$f^{\underline{k}}(0)=0$ for any~$k$ odd (resp., even). For instance, Bingham alternatives associated with~$f(s)=\exp(s^2)$, for which $f^{\underline{k}}(0)=0$ for any~$k$ odd, cannot be detected at any polynomial rate by finite Sobolev tests whose non-zero coefficients~$v_k$ show at odd ranks~$k$ only (such as, e.g., the Rayleigh test).

Finally, we point out a further interesting practical consequence of this parity feature in Theorem~\ref{Sobolfinbisbis}. Any finite Sobolev test for which both $v_1$ and $v_2$ are non-zero (such as, e.g., the Sobolev test whose statistic is the sum of the Rayleigh and Bingham statistics) will exhibit the same consistency rate as the oracle Sobolev test based on $v_k = \delta_{k,\ell}$, where $\ell$ is the smallest positive integer for which~$f^{\underline{\ell}}(0)\neq 0$. Therefore, a Sobolev test for which both $v_1$ and $v_2$ are non-zero will always be rate-consistent against the considered local~alternatives.

%-------------------------------%
\subsection{Infinite Sobolev tests}
%-------------------------------%

%
We now turn our attention to infinite Sobolev tests, based on statistics of the form
\begin{align}
\label{untruncSobol}
S^{(n)}_v
:=
\sum_{k=1}^{\infty} v_k^2 \bigg\| \frac{1}{\sqrt{n}} \sum_{i=1}^n {\bf G}_{p,k}(\Ub_i) \bigg\|^2
\end{align}
(see~\eqref{Sobolgen}), where infinitely many coefficients~$v_k$ are non-zero. Recall that these tests reject the null hypothesis of uniformity at asymptotic level~$\alpha$ when~$S^{(n)}_v$ exceeds the upper $\alpha$-quantile of the distribution of~$\sum_{k=1}^{\infty} v_k^2 \chi^2_{d_{p,k}}$. %
If~$\sum_{k=1}^{\infty}
v_k^2
d_{p,k}
< \infty
$, then these tests indeed lead to rejection with asymptotic probability~$\alpha$ under the null hypothesis; see Theorem~\ref{asymptnull}.

Now, as mentioned in the introduction, infinite Sobolev tests are consistent under any fixed alternative to uniformity as soon as~$v_k\neq 0$ for any~$k$; see Theorem~4.4 in \cite{Gine1975}. Yet asymptotic local powers, or even simply consistency rates, are not precisely available in the literature, which calls for an extension of the results of the previous section to such infinite Sobolev tests. We have the following result.

\begin{theorem}
\label{Sobolfininfin}
Fix an integer~$p\geq 2$ and a sequence~$v=(v_k)$ satisfying %
$ \sum_{k=1}^{\infty}
v_k^2
d_{p,k}
< \infty$. Let~$f$ be an angular function that is infinitely many times differentiable at zero. Let~$k_v:=\min\{k:v_k\neq 0\}$, fix~$\tau>0$, and let~$\xi_{p,k}(\tau)$ be as in~\eqref{ksith41}.
Then, %
\begin{enumerate}[label=(\textit{\roman{*}}),ref=\textit{\roman{*}}]
	\item under~${\rm P}\n_{\kappa_n, f}$, with $\kappa_n=o(n^{-1/(2{k_v})})$,
	$
	S^{(n)}_v
	\stackrel{\mathcal{D}}{\to}
	\sum_{k=1}^{\infty} v_k^2 \chi^2_{d_{p,k}}
	$; \label{Sobolfininfin:1}
	\item under~${\rm P}\n_{\kappa_n, f}$, with $\kappa_n=n^{-1/(2{k_v})} \tau$,
	$
	S^{(n)}_v
	\stackrel{\mathcal{D}}{\to}
	\sum_{k=1}^{\infty} v_k^2 \chi^2_{d_{p,k}}(\delta_{k,k_v}\xi_{p,k_v}(\tau))
	$;\label{Sobolfininfin:2}
	\item under~${\rm P}\n_{\kappa_n, f}$, where~$n^{1/(2{k_v})}\kappa_n\to\infty$ and~$f^{\underline{k_v}}(0)\neq 0$,
	$
	{\rm P}\n_{\kappa_n,f}[S^{(n)}_v>c]\to 1
	$
	for any~$c>0$, %
	so that the corresponding Sobolev test is consistent under this sequence of alternatives.\label{Sobolfininfin:3}
\end{enumerate}
\end{theorem}

This result is the analogue of Theorem~\ref{Sobolfin} for Sobolev tests based on an infinite collection of spherical harmonics. Clearly, it states that, if %
$\sum_{k=1}^{\infty}
v_k^2
d_{p,k}
< \infty
$ holds (note that this condition is also needed under the null hypothesis) and if~$f$ is infinitely many times differentiable at zero, then the detection threshold of an infinite Sobolev test is~$\kappa_n\sim n^{-1/(2k_v)}$, provided that~$f^{\underline{k_v}}(0)\neq 0$. Inspection of the proof of Theorem~\ref{Sobolfininfin} reveals that Theorem~\ref{Sobolfinbisbis} (hence also Corollary~\ref{Sobolfinbis}) can be similarly adapted to infinite Sobolev tests, which will allow one to determine the detection threshold also in the case where~$f^{\underline{k_v}}(0)=0$.

%

%-------------------------------%
\section{Numerical illustrations}
\label{sec:simus}
%-------------------------------%

The objective of this section is to conduct Monte Carlo exercises in order to illustrate the various results obtained in the paper. We focus on the problem of testing for uniformity obtained with~$\Gamb={\bf I}_d$ and consider three types of alternatives that allow us to cover the subtle non-null results described in the previous sections.

\begin{figure}[!htb]
%\vspace*{-0.5cm}
\centering
\includegraphics[width=0.85\textwidth]{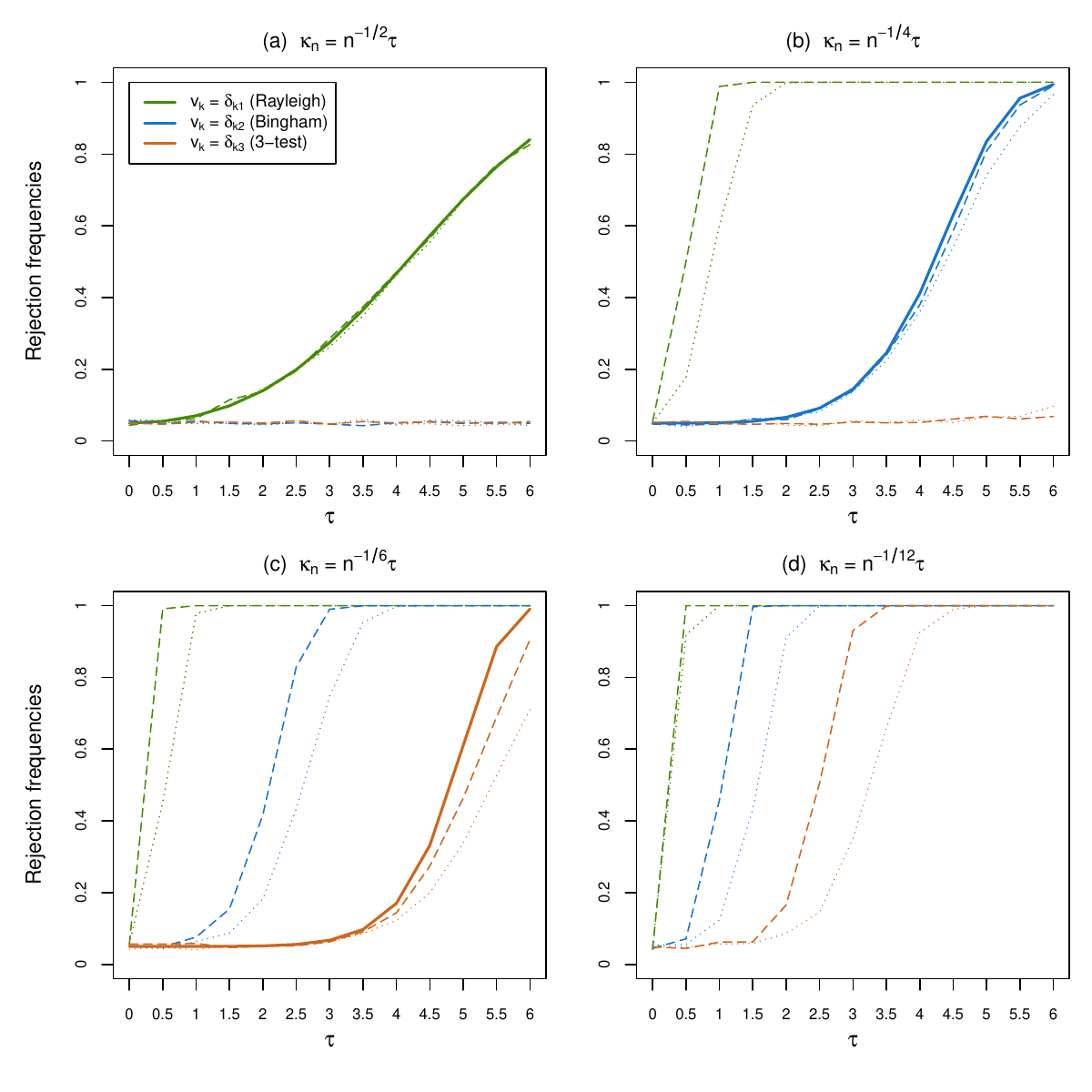}
\vspace*{-0.25cm}
\caption{\small Rejection frequencies of the Rayleigh test (green), Bingham test (blue), and $3$-test (dark orange), all conducted at asymptotic level~$\alpha=5\%$ for samples generated by the rotationally symmetric model with angular function~$f(s)=\exp(s)$ and concentration~$\kappa_n=n^{-1/\ell}\tau$, for $\ell=2,4,6,12$ and~$\tau\in\{0,0.5,1,\ldots,6\}$. %
	Whenever the rejection frequencies are non-trivial (that is, strictly between~$\alpha$ and~$1$), the corresponding asymptotic powers are also provided in solid lines.}
\label{Fig1}
\end{figure}

First, for each sample size~$n\in\{500,2,\!500\}$, each~$\ell\in\{2,4,6,12\}$, and each~$\tau\in \{0,0.5,1,\ldots,6\}$, we generated~$M=5,\!000$ mutually independent random samples $\Ub_1, \ldots, \Ub_n$ from the von Mises--Fisher distribution ($f(s)=\exp(s)$) with location~$\thetab_0=(0,0,1)'\in\mathcal{S}^2$ and concentration~$\kappa_n=n^{-1/\ell}\tau$. In each of the resulting samples, we performed the following tests at asymptotic level~$\alpha=5\%$: the Rayleigh test of uniformity ($\Gamb={\bf I}_3$), the Bingham test of uniformity, and the ``$3$-test'' of uniformity. These are the Sobolev tests associated with~$v_k=\delta_{k,1}$, $v_k=\delta_{k,2}$, and~$v_k=\delta_{k,3}$, respectively. %
The resulting empirical rejection frequencies are plotted as functions of~$\tau$ in Figure~\ref{Fig1}, which also provides the corresponding asymptotic power curves obtained from Theorem~\ref{Sobolfin} whenever these are non-trivial (that is, are strictly between~$\alpha$ and~$1$).
In the present von Mises--Fisher case, $f^{\underline{k}}(0)\neq 0$ for any~$k$, so that the Rayleigh test, the Bingham test, and the $3$-test have a detection threshold given by~$\kappa_n \sim n^{-1/2}$, $\kappa_n \sim n^{-1/4}$, and~$\kappa_n \sim n^{-1/6}$, respectively (Theorem~\ref{Sobolfin}). This is perfectly supported by Figure~\ref{Fig1}, where the Rayleigh test first shows power in Panel~(a), whereas the Bingham test and the $3$-test first show power in Panels~(b) and~(c), respectively. For alternatives that are over their respective detection threshold, the various tests exhibit rejection frequencies that are compatible, at the finite-sample sizes considered, with consistency. Finally, at their respective detection threshold, the agreement between rejection frequencies and asymptotic powers is excellent (but maybe for the $3$-test in Panel~(c), yet, as we checked by performing further simulations, this improves for larger sample sizes).

\begin{figure}[!htb]
\centering
\includegraphics[width=0.85\textwidth]{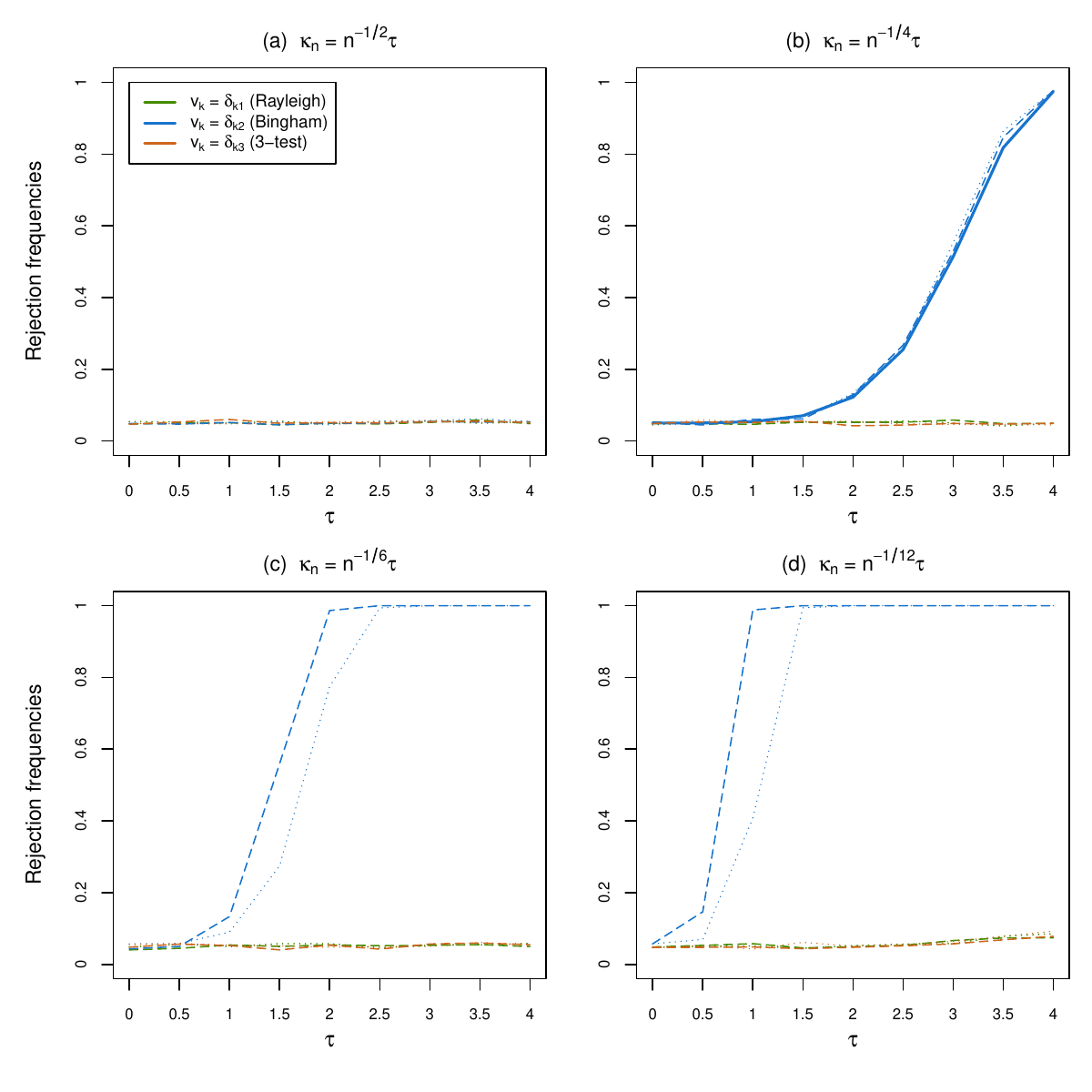}
\vspace*{-0.25cm}
\caption{\small Same experiment as in Figure~\ref{Fig1}, but with samples generated from a rotationally symmetric model with angular function~$f(s)=\exp(s^2)$ and~$\tau\in\{0,0.5,1,\ldots,4\}$.}
\label{Fig2}
\end{figure}

Second, we repeated the same exercise as above but with Watson distributions replacing the von Mises--Fisher ones, that is, with~$f(s)=\exp(s^2)$ instead of~$f(s)=\exp(s)$. In this exercise, we rather used~$\tau\in \{0,0.5,1,\ldots,4\}$; this discrepancy with the values of~$\tau$ in the first numerical exercise above is unimportant (the maximal value of~$\tau$ was chosen in each simulation to obtain non-trivial asymptotic powers that roughly end up at one at the detection threshold; see Panel~(b) of Figure~\ref{Fig2}). In this Watson case, $f^{\underline{k}}(0) \neq 0$ if and only if~$k$ is even; thus, the detection threshold of the Bingham test is still~$\kappa_n \sim n^{-1/4}$ (Theorem~\ref{Sobolfin}), while the Rayleigh test and the $3$-test should be blind to all alternatives considered (Corollary~\ref{Sobolfinbis}). This is perfectly in line with the empirical rejection frequencies plotted in Figure~\ref{Fig2}.

\begin{figure}[!htb]
%\vspace*{-0.5cm}
\centering
\includegraphics[width=0.85\textwidth]{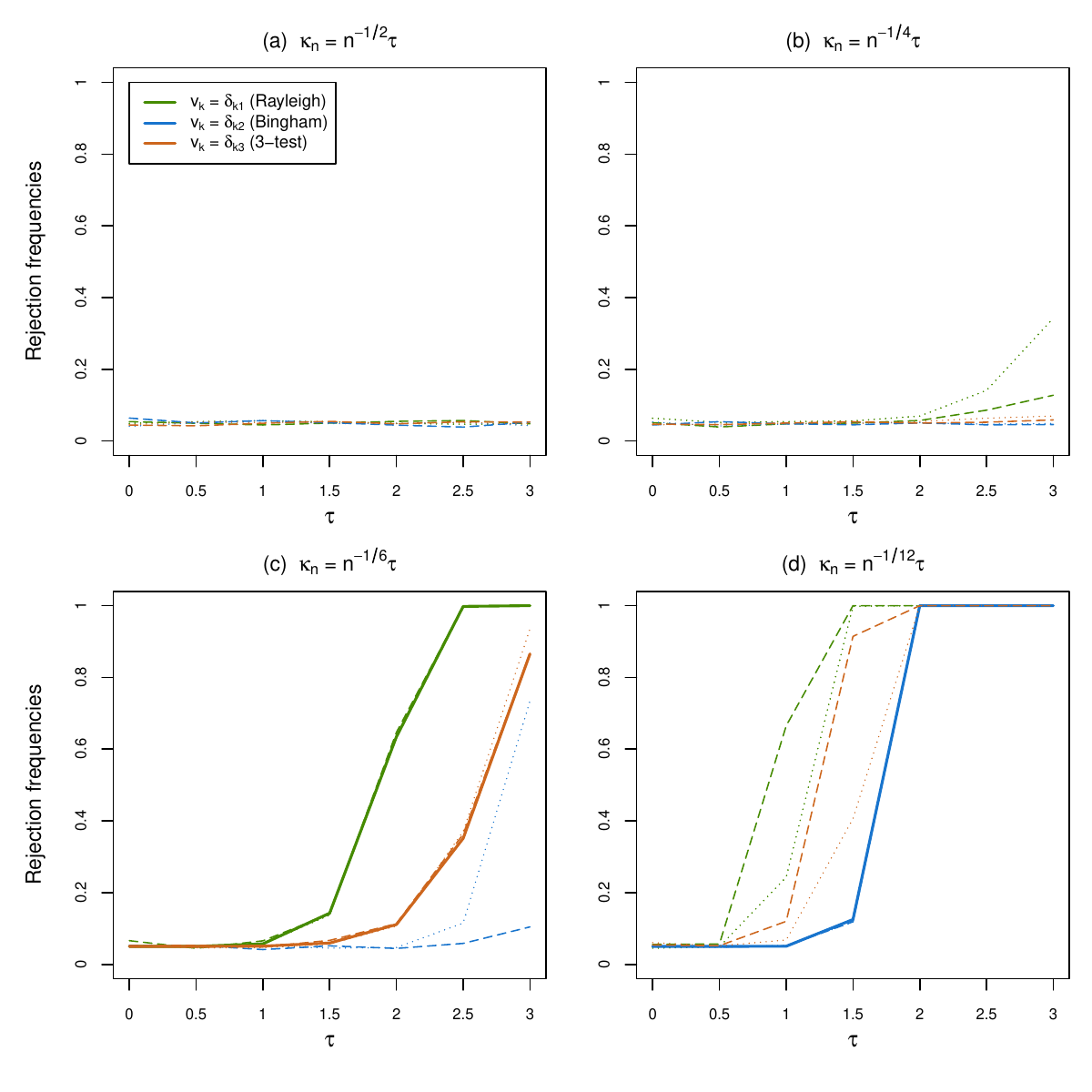}
\vspace*{-0.25cm}
\caption{\small Same experiment as in Figure~\ref{Fig1}, but with samples generated from a rotationally symmetric model with angular function~$f(s)=\exp(s^3)$ and~$\tau\in\{0,0.5,1,\ldots,3\}$.}
\label{Fig3}
\end{figure}

Third, we conducted a last simulation where observations were rather generated using the angular function~$f(s)=\exp(s^3)$ (and~$\tau\in \{0,0.5,1,\ldots,3\}$). The resulting empirical rejection frequencies are plotted in Figure~\ref{Fig3}. Here, $f^{\underline{k}}(0) \neq 0$ if and only if~$k$ is a multiple of~$3$. The detection threshold of the $3$-test and of the Rayleigh test should then be~$\kappa_n\sim n^{-1/6}$ (Theorems~\ref{Sobolfin} and~\ref{Sobolfinbisbis}, respectively), whereas the one of the Bingham test should be~$\kappa_n\sim n^{-1/12}$ (Theorem~\ref{Sobolfinbisbis}). Clearly, this is confirmed by the rejection frequencies in Figure~\ref{Fig3}, where the agreement with the corresponding asymptotic powers (see Panel~(c) for the Rayleigh test and the $3$-test, and Panel~(d) for the Bingham test) is, again, excellent.

We conclude that all numerical experiments above (\textit{i}) strongly support our asymptotic results
and (\textit{ii}) further reveal that the asymptotic behavior most often kicks in at sample sizes as small as~$n=500$, which, in view of the highly non-standard nature of the results (e.g., detection thresholds in~$\kappa_n\sim n^{-1/12}$), is quite remarkable.

%-------------------------------%
\section{Symmetry of comet orbits}
\label{sec:data}
%-------------------------------%

%
Orbits of celestial bodies, such as planets and comets, have attracted scientists' attention for a long time. \cite{Bernoulli1735} already discussed whether the clustering of the planets' orbits about the ecliptic, nowadays explained by their origin in the protoplanetary disk, could have happened ``by chance''. The study of comet orbits has been more intricate. Long-period comets (with periods larger than $200$ years) are thought to arise from the roughly spherical Oort cloud, containing icy planetesimals that were ejected from protoplanetary disks by giant planets. These icy planetesimals became heliocentric comets when their orbits were affected by random perturbations of passing stars and the galactic tide (see, e.g., Sections 5 and 7.2 in \citet{Dones2015} and references therein). This conjectured past of the Oort cloud explains the nearly isotropic distribution of long-period comets \citep[evidenced, e.g., in][]{Wiegert1999}, sharply contrasting with the ecliptic-clustered\nopagebreak[4] orbits of short-period comets originating in the flattened Kuiper belt (see Figure \ref{fig:kuiper}).

As illustrated in \cite{Watson1970} and 
\cite{Jupp2003}, assessing the uniformity of orbits can be formalized as testing the uniformity on $\mathcal{S}^2$ of their directed unit normal vectors. An orbit with \textit{inclination} $i\in[0,\pi]$ and \textit{longitude of the ascending node} $\Omega\in[0,2\pi)$ \cite[see][]{Jupp2003} has directed normal vector $(\sin(i)\sin(\Omega),-\sin(i)\cos(\Omega),\cos(i))'$ to the orbit's plane (see the illustrative graphs in Figure \ref{fig:comets}). The sign of the vector reflects if the orbit is prograde (the normal vector points to the northern hemisphere) or retrograde (southern hemisphere). Following the data acquisition methodology from \cite{Cuesta-Albertos2009}, we retrieved the long-period elliptic-type single-apparition comets present in the JPL Small-Body Database Search Engine\footnote{\url{https://ssd.jpl.nasa.gov/sbdb_query.cgi}}. As of 2022-05-28, $n=610$ comets were retrieved, excluding the comet fragments that create spurious clusters. The dynamic nature of the database, with additions of first-ever observed comets and updates on the data for former comets, generated the noticeable difference between the $208$ comets analyzed in \cite{Cuesta-Albertos2009}. The same methodology was used to retrieve a dataset of $n=784$ short-period comets. Both datasets are contained in the \texttt{comets} object in the \texttt{sphunif} R package \citep{Garcia-Portugues-sphunif}.

\begin{figure}[h!]
\centering
\begin{subfigure}{0.5\textwidth}
	\centering
	\includegraphics[width=0.5\textwidth,trim={3cm 0cm 3cm 0cm},clip]{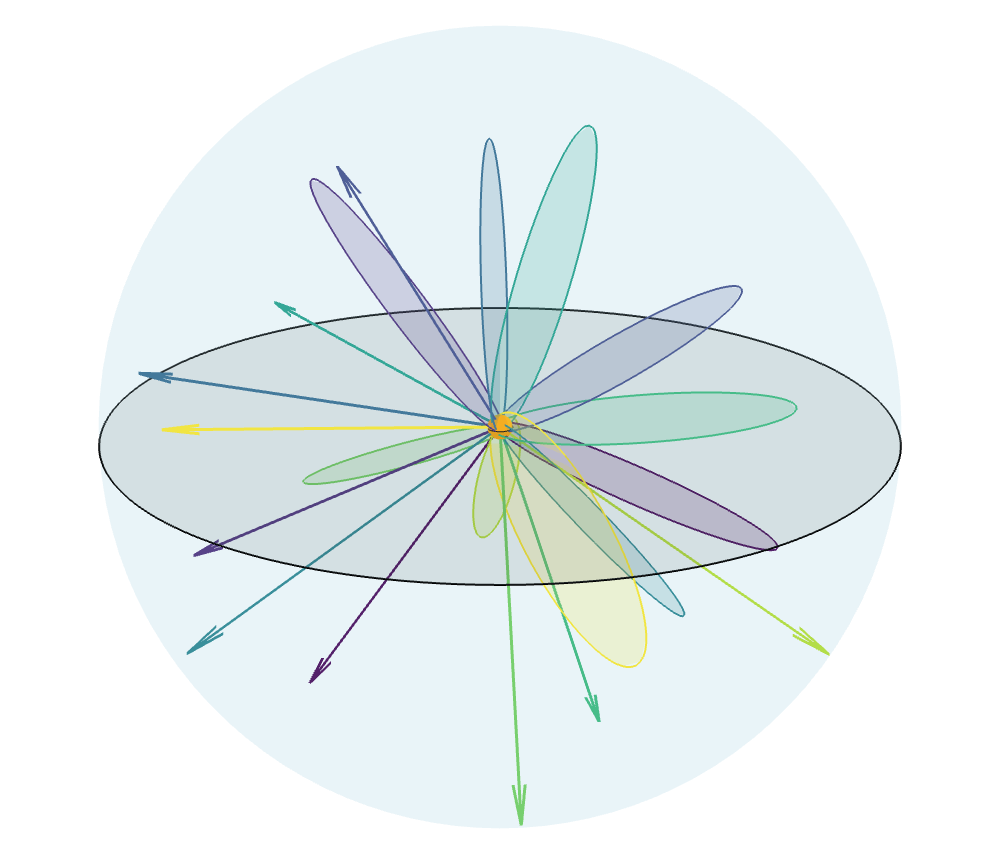}\includegraphics[width=0.5\textwidth,trim={3cm 0cm 3cm 0cm},clip]{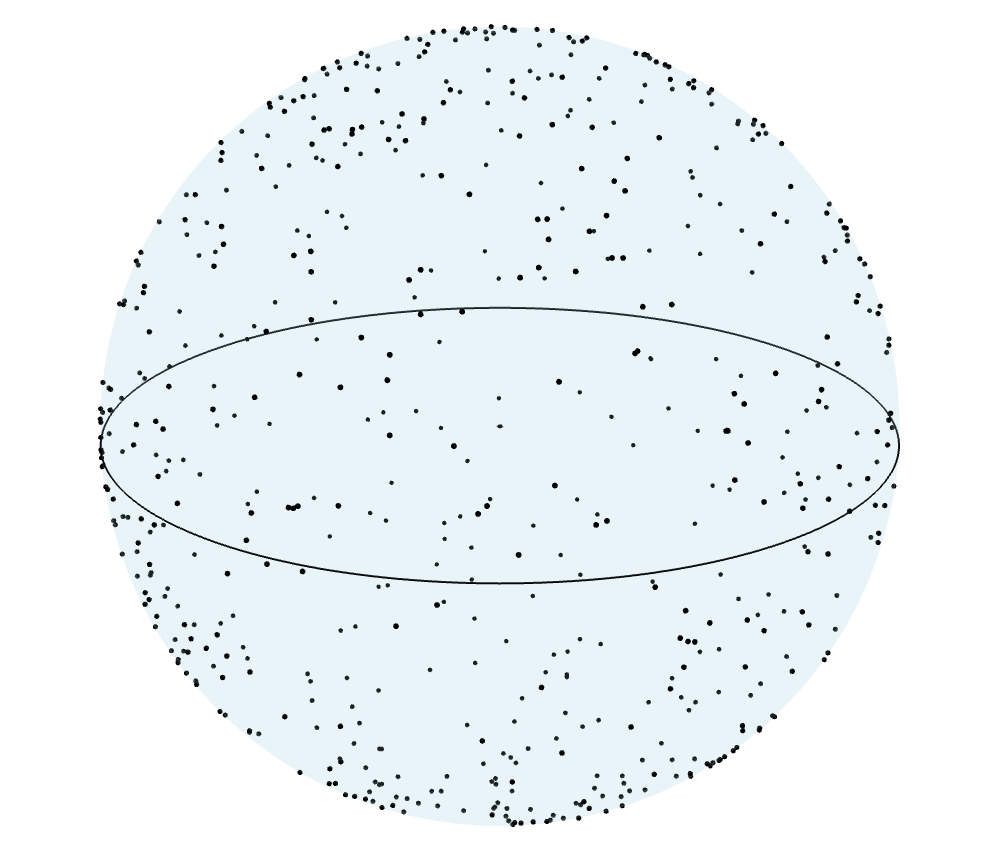}
	\caption{Long-period comets} 
	\label{fig:oort}
\end{subfigure}%
\begin{subfigure}{0.5\textwidth}
	\includegraphics[width=0.5\textwidth,trim={3cm 0cm 3cm 0cm},clip]{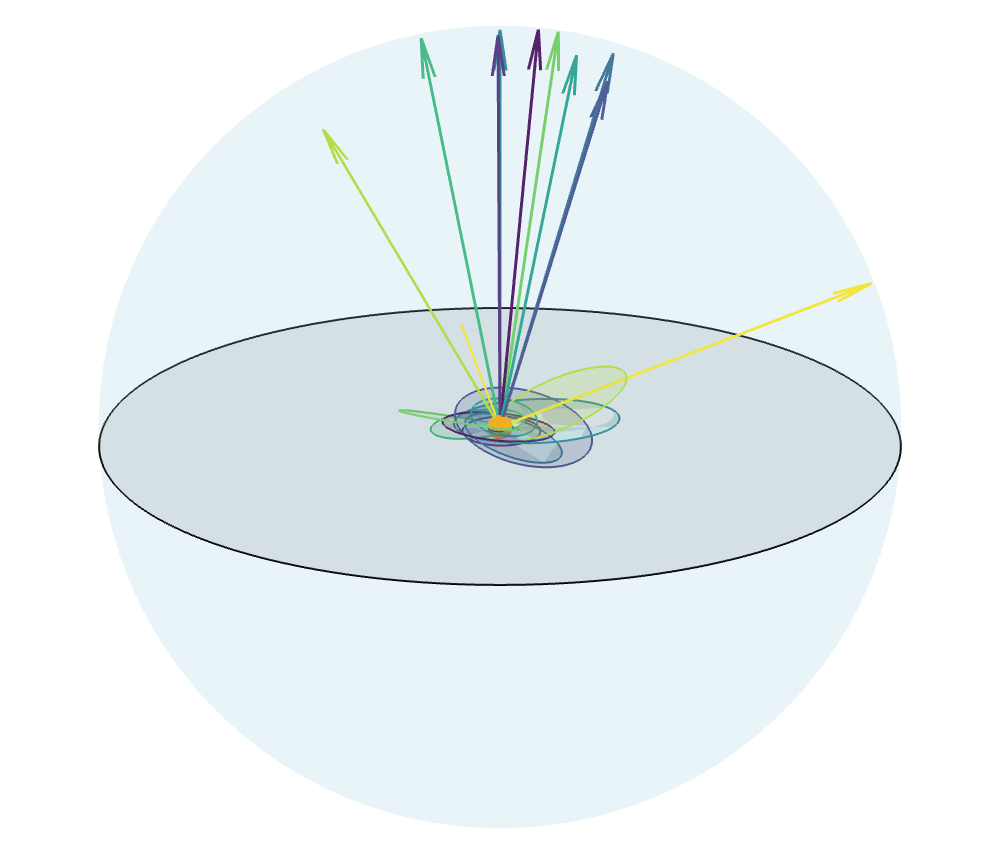}\includegraphics[width=0.5\textwidth,trim={3cm 0cm 3cm 0cm},clip]{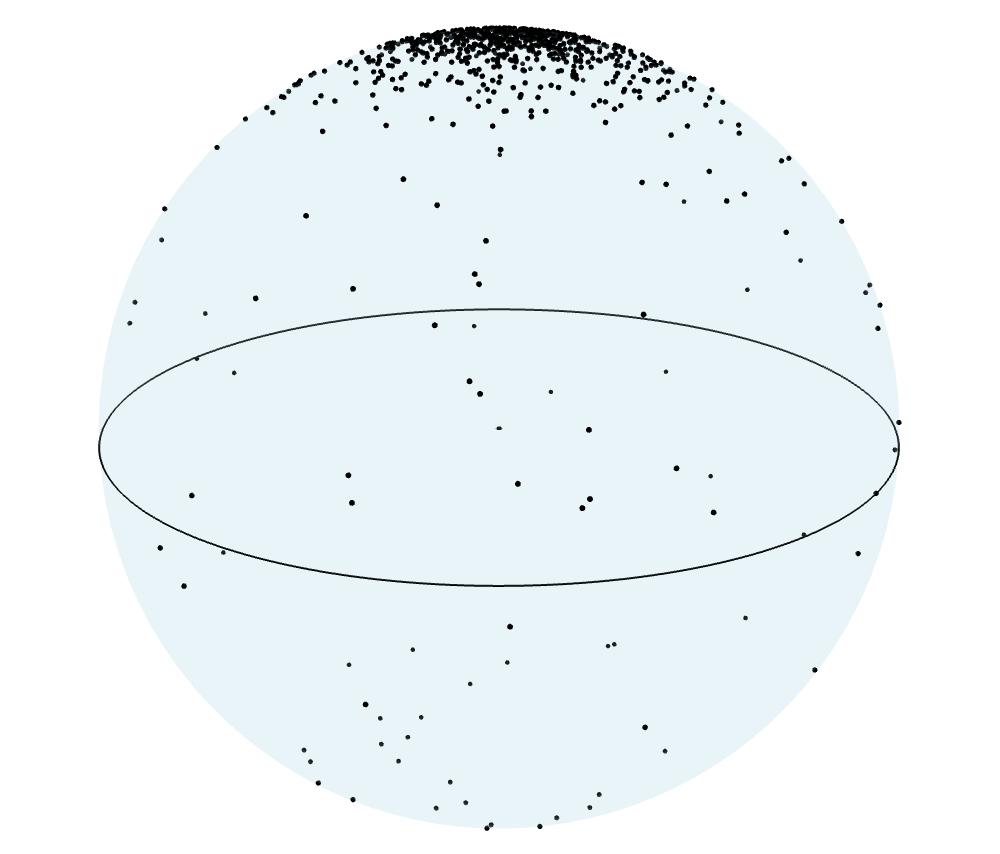}
	\caption{Short period comets}
	\label{fig:kuiper}
\end{subfigure}%
\caption{\small Orbits of long- and short-period comets and their normal vectors. Within each figure, the left plot displays ten illustrative elliptical orbits and their associated normal vectors, with the ecliptic plane shown in gray and the Sun represented as an orange sphere (one focus of the elliptical orbits). The right plot in each figure shows the full dataset of normal orbit vectors on~$\mathcal{S}^2$.}
\label{fig:comets}
\end{figure}

We carried out finite Sobolev tests with weights $v_k=\delta_{k,k_v}$ for $k_v=1,\ldots,6$ to these the uniformity of orbits of long-period comets. To complement the outcome of these tests, we derived ``indicative powers'' using Theorem~\ref{Sobolfin}\eqref{Sobolfin:2}. For a Sobolev test based on $S_v^{(n)}$, these indicative powers were computed by: (a) considering the exponential-like alternatives ${\bf x}\mapsto c_{p,\kappa,b}\exp(\kappa({\bf x}'\boldsymbol\mu)^b)$ for $b=1,\ldots,6$; (b) computing maximum likelihood estimators $(\hat{\boldsymbol{\mu}}_b,\hat\kappa_b)$; (c) fitting the sequence of alternatives ${\rm P}\n_{\kappa_n, f_{\kappa,b}}$ to the sample by setting $\hat\tau_b:=n^{1/(2k_v)}\hat{\kappa}_b$; (d) obtaining the rejection probability under the fitted alternative as $\mathrm{P}\big[\sum_{k=1}^{K_v} v_k^2 \chi^2_{d_{p,k}}(\delta_{k,k_v}\xi_{p,k_v}(\hat\tau_b))>c_\alpha\big\vert\hat\tau_b\big]$, where $c_\alpha$ is the upper $\alpha$-quantile of the distribution of~$\sum_{k=1}^{K_v} v_k^2 \chi^2_{d_{p,k}}$. Table \ref{tab:comets:1} reports the $p$-values of the tests and the maximum of their indicative powers under the alternatives with $b=1,\ldots,6$ at the asymptotic level $\alpha=5\%$.

\begin{table}[h!]
\centering
\begin{tabular}{c|c|cccccc}
	\toprule
	Symmetry hypothesis & Probability & $k_v=1$ & $k_v=2$ & $k_v=3$ & $k_v=4$ & $k_v=5$ & $k_v=6$\\\midrule
	Uniformity & $p$-value & $0.2120$ & $0.0000$ & $0.0793$ & $0.0015$ & $0.5904$ & $0.3821$\\
	& Power & $0.4013$ & $0.8510$ & $0.0504$ & $0.0762$ & $0.0500$ & $0.0519$\\ \midrule
	Rotational symmetry & $p$-value & $0.1685$ & $0.5820$ & $0.1019$ & $0.2797$ & $0.6787$ & $0.7185$\\
	about $\boldsymbol{\mu}_1=(0,0,1)'$ & Power & $0.3757$ & $0.0512$ & $0.0500$ & $0.0500$ & $0.0500$ & $0.0500$\\ \bottomrule
\end{tabular}
\caption{\small Asymptotic $p$-values and maximum indicative powers for the $S_v^{(n)}$-based tests of symmetry of directions applied to the normal vectors of the orbits of long-period comets. The tests were conducted at the asymptotic level $\alpha=5\%$ for the sequence of weights $v_k=\delta_{k,k_v}$.}
\label{tab:comets:1}
\end{table}

Significant non-uniformity is detected in the orbits of long-period comets for the spherical harmonics associated with $k_v=2$ (Bingham test; highest power) and $k_v=4$ (marginal power). Along with the non-rejections for odd $k_v=1,3,5$ (with $\alpha=5\%$), the outcome of the analysis signals a non-uniformity in the form of a low-order axial alternative. The observational bias of long-period comets, as described in \cite{Jupp2003}, may explain this rejection: search programs prioritize observational bandwidths near the ecliptic plane to maximize the discovery of (short- or long-period) comets. Near-ecliptic-plane orbits correspond to symmetric clusters of normal vectors on the north and south of $\mathcal{S}^2$, which are well detected by even spherical harmonics.

Despite uniformity being rejected, the distribution of long-period orbits is markedly symmetric. Indeed, Table \ref{tab:comets:1} shows that no significant deviation from rotational symmetry about $\boldsymbol{\mu}_1=(0,0,1)'$, the normal to the ecliptic plane, is observed (this kind of symmetry is compatible with the presence of uniform observational bias near the ecliptic). However, rotational symmetry about $\boldsymbol{\mu}_1=(0,0,1)'$ is emphatically rejected for short-period comets by the Rayleigh test, as Table \ref{tab:comets:2} reveals. A closer examination of the rejection cause reveals that this is attributed to a minor shift of the rotational symmetry axis. If running the unspecified-$\boldsymbol{\mu}$ variants of the Rayleigh and Bingham tests from \cite{Garcia-Portugues2020}, namely $\phi_\mathrm{vMF}^\mathrm{loc}$ and $\phi_\dagger^\mathrm{sc}$, rotational symmetry is not rejected, with a resulting estimated rotational symmetry direction $\hat{\boldsymbol{\mu}}_1=(0.0355,-0.0144,0.9993)\pr$.

\begin{table}[h!]
\centering
\begin{tabular}{c|cccccc}
	\toprule
	Rotational symmetry & $k_v=1$ & $k_v=2$ & $k_v=3$ & $k_v=4$ & $k_v=5$ & $k_v=6$\\\midrule
	About $\boldsymbol{\mu}_1=(0,0,1)'$ & $0.0006$ & $0.7957$ & $0.8533$ & $0.0317$ & $0.6338$ & $0.1607$\\
	About $\boldsymbol{\mu}_1$ unspecified & $0.2092$ & $0.9229$ & --- & --- & --- & --- \\ \bottomrule
\end{tabular}
\caption{\small Asymptotic $p$-values for the $S_v^{(n)}$-based tests of rotational symmetry applied to the normal vectors of the orbits of short-period comets. In the second row, the unspecified-$\boldsymbol{\mu}_1$ tests $\phi_\mathrm{vMF}^\mathrm{loc}$ ($k_v=1$) and $\phi_\dagger^\mathrm{sc}$ ($k_v=2$) from \cite{Garcia-Portugues2020} are used. The tests were conducted at the asymptotic level $\alpha=5\%$ for the sequence of weights $v_k=\delta_{k,k_v}$.}
\label{tab:comets:2}
\end{table}

%-------------------------------%
\section{Perspectives for future research}
\label{sec:perspectives}
%-------------------------------%

The present paper studied the asymptotic powers of a class of Sobolev tests under natural alternatives. The results require minimal assumptions (only differentiability assumptions of the angular function~$f$ at zero are imposed) and they cover the whole class of Sobolev tests. The consistency rate of an arbitrary Sobolev test and the corresponding asymptotic powers were obtained under \emph{almost all} rotationally symmetric distributions. To be more specific, if the angular function~$f$ is infinitely many times differentiable at zero with~$f^{\underline{k}}(0)=0$ for any positive integer~$k$ (consider for instance the function~$f$ such that~$f(s)=1+\exp(-1/s^2)$ for~$s\neq 0$ and~$f(0)=1$), then our results only state that all Sobolev tests will be blind to local alternatives in~$\kappa_n\sim n^{-1/(2q)}$ irrespective of the positive integer~$q$. As a corollary, consistency rates in this particular case remain unknown (consistency follows from Theorem~4.4 in \citealp{Gine1975}). The consistency rate and resulting local powers could be explored in future work, even though this may be of academic interest only (one may argue that if a test does not show polynomial consistency rates in a specific parametric model, then this test is in any case of low practical relevance in that model).

Another avenue for future research is to consider other distributional setups. Rotationally symmetric distributions involve a single particular direction, namely the direction of~$\thetab$, hence are ``single-spiked'' distributions. More general distributions in directional statistics are rather ``multi-spiked'' ones. For instance, the \cite{Bingham1974} distributions are multi-spiked extensions of the rotationally symmetric distributions associated with~$f(s)=\exp(s^2)$, that is, extensions of the Watson distributions. It would be of interest to extend the results of the present work by characterizing the local non-null behavior of Sobolev tests against such multi-spiked alternatives. This is by no means easy since, in principle, each ``spike'' may involve its own concentration parameter~$\kappa$. Another possible direction for future research is to go away from the fixed-$p$ setup considered in this paper to consider high-dimensional asymptotic scenarios under which the dimension~$p=p_n$ diverges to infinity as~$n$ does. Consistency rates and local powers in this case are available for very specific Sobolev tests only; we refer to \cite{Cutting2017,Cutting2021}, for the Rayleigh test and for the Bingham test, respectively. Results that would apply to arbitrary Sobolev tests in high dimensions would be of interest. This will be considered in future work.

%-----------------------------------------------%
\section*{Supplementary materials}
%-----------------------------------------------%

Supplementary materials provide the proofs of the stated results.

%-----------------------------------------------%
\section*{Acknowledgments}
%-----------------------------------------------%

The first author acknowledges support by grant PID2021-124051NB-I00, funded by MCIN/\-AEI/\-10.13039/\-501100011033 and by ``ERDF A way of making Europe''. The first author acknowledges support from ``Convocatoria de la Universidad Carlos III de Madrid de Ayudas para la recualificación del sistema universitario español para 2021--2023'', funded by Spain's Ministerio de Ciencia, Innovación y Universidades. The second author's research is supported by a research fellowship from the Francqui Foundation and by the Program of Concerted Research Actions (ARC) of the Universit\'{e} libre de Bruxelles. The third author's research is supported by the ARC Program of the Universit\'{e} libre de Bruxelles and by the Fonds Thelam from the Fondation Roi Baudouin.

%\bibliographystyle{apalike}
%\bibliography{DirectionalStats.bib,no-DirectionalStats.bib}
% DirectionalStats.bib v1.25 (2021-07-20)

\fi

\ifsupplement

\newpage
%-----------------------------------------------%
\title{Supplementary materials for ``On a class of Sobolev tests for symmetry of directions, their detection thresholds, \newline and asymptotic powers''}
\setlength{\droptitle}{-1cm}
\predate{}%
\postdate{}%
\date{}
%-----------------------------------------------%

%-----------------------------------------------%
\author{Eduardo Garc\'ia-Portugu\'es$^{1}$, Davy Paindaveine$^{2,3}$, and Thomas Verdebout$^{2}$}
\footnotetext[1]{Department of Statistics, Universidad Carlos III de Madrid (Spain).}
\footnotetext[2]{D\'{e}partement de Math\'{e}matique and ECARES, Universit\'{e} libre de Bruxelles (Belgium).}
\footnotetext[3]{Corresponding author. e-mail: \href{mailto:dpaindav@ulb.ac.be}{dpaindav@ulb.ac.be}.}
\maketitle
%-----------------------------------------------%

%-----------------------------------------------%
\begin{abstract}
	These supplementary materials provide the proofs of the stated results in the paper.
\end{abstract}
\begin{flushleft}
	\small\textbf{Keywords:} Consistency rates; detection thresholds; directional statistics; hypothesis tests.
\end{flushleft}
%-----------------------------------------------%

\appendix

%-------------------------------%
\section{Proofs of Section \ref{sec:moments}}
%-------------------------------%

The proof of Proposition \ref{Exprop} requires the following result, that is classical in spherical harmonic analysis; see, e.g., Theorem~1.2.9 in \cite{Dai2013|SM}.

\begin{thm}[Funk--Hecke Theorem]\label{thm:funk}
	Let~$p\geq 3$ be an integer and~$k$ be a nonnegative integer. Let $h:[-1,1]\to\mathbb{R}$ be a measurable function such that~$	\int_{-1}^1 (1-s^2)^{(p-3)/2} |h(s)| \,\mathrm{d}s$ is finite. Then, for any~$\varphi_k \in\mathcal{H}^p_k$ and any $\etab\in \mathcal{S}^{p-1}$,
	\begin{align*}
		\int_{\mathcal{S}^{p-1}} h(\xib'\etab) \varphi_k(\xib) \,\mathrm{d}\sigma_p(\xib)=\lambda_k \varphi_k(\etab),
	\end{align*}
	with
	\begin{align*}
		\lambda_k=
		\frac{\omega_{p-2}}{C_k^{(p-2)/2}(1)}
		\int_{-1}^1 (1-s^2)^{(p-3)/2} h(s) C_k^{(p-2)/2}(s) \,\mathrm{d}s,
	\end{align*}
	where~$\omega_{p-2}$ is the surface area of~$\mathcal{S}^{p-2}$.
\end{thm}

\begin{proof}[Proof of Proposition \ref{Exprop}]
	Throughout, we write~$\Ub:=\Ub_1$. Let us first prove~\eqref{shiftstep1} for~$p=2$, where $\Ub=(U_1,U_2)'$. We have
	\begin{align*}
		{\rm E}\n_{\kappa_n,f}[g_{1,k}(\Ub)]
		&=
		\sqrt{2} {\rm E}\n_{\kappa_n,f}[
		\cos(k \arccos U_2) \mathbb{I}[U_1\geq 0]
		+
		\cos(k (2\pi-\arccos U_2)) \mathbb{I}[U_1<0]
		]
		\\
		&
		=
		\sqrt{2}
		{\rm E}\n_{\kappa_n,f}[
		\cos(k \arccos U_2)]
		=
		\sqrt{2}
		{\rm E}\n_{\kappa_n,f}[C_k^{0}(U_2)]
		=
		\sqrt{2}
		{\rm E}\n_{\kappa_n,f}[C_k^{0}(\Ub'\thetab_0)]
		.
	\end{align*}
	Similarly, since~$U_1\neq 0$ almost surely under~${\rm P}\n_{\kappa_n,f}$,
	\begin{align*}
		{\rm E}\n_{\kappa_n,f}[g_{2,k}(\Ub)]
		&=
		\sqrt{2} {\rm E}\n_{\kappa_n,f}[
		\sin(k \arccos U_2) \mathbb{I}[U_1> 0]
		+
		\sin(k (2\pi-\arccos U_2)) \mathbb{I}[U_1<0]
		]
		\\
		&=
		\sqrt{2} {\rm E}\n_{\kappa_n,f}[
		\sin(k \arccos U_2) (\mathbb{I}[U_1> 0]-\mathbb{I}[U_1<0])
		]
		=
		0
		,
	\end{align*}
	where the last equality results from the fact, the distribution being rotationally symmetric about~$\thetab_0=(0,1)'$, we have that~$(-U_1,U_2)$ and~$(U_1,U_2)$ are equal in distribution under~${\rm P}\n_{\kappa_n,f}$.
	This establishes~\eqref{shiftstep1} for~$p=2$. For~$p\geq 3$, taking $h(s)=c_{p,\kappa_n,f} f(\kappa_n s)$ in Theorem \ref{thm:funk} yields, for $r=1, \ldots, d_{p,k}$,
	\begin{align*}
		\mathrm{E}\n_{\kappa_n,f}[{g}_{r,k}(\Ub)]
		&=
		\frac{1}{\omega_{p-2}}
		\int_{\mathcal{S}^{p-1}}
		h(\xb\pr \thetab_0){g}_{r,k}(\xb)
		\,\mathrm{d}\sigma_p(\xb)
		=
		\frac{\lambda_k {g}_{r,k}(\thetab_0)}{\omega_{p-2}}
		\\
		&=
		\frac{1}{C_k^{(p-2)/2}(1)}
		\bigg(
		c_{p,\kappa_n,f}
		\int_{-1}^1 (1-s^2)^{(p-3)/2} f(\kappa_n s) C_k^{(p-2)/2}(s) \,\mathrm{d}s
		\bigg)
		{g}_{r,k}(\thetab_0)
		\\
		&=
		\frac{1}{C_k^{(p-2)/2}(1)}
		\mathrm{E}\n_{\kappa_n,f}[C_k^{(p-2)/2}(\Ub'\thetab_0)]
		{g}_{r,k}(\thetab_0)
		.
	\end{align*}
	Recalling that~$\thetab_0=(0,\ldots,0,1)'\in\R^p$,~\eqref{shiftstep1} then follows from the fact that
	\begin{align*}
		g_{1,k}(\thetab_0)
		&=
		\frac{1}{\sqrt{d_{p,k}}}
		\Big(1+ \frac{2k}{p-2}\Big)
		C_k^{(p-2)/2}(1),
		\\
		g_{r,k}(\thetab_0)&=0
		,
		\quad
		r=2,\ldots,d_{p,k}
		;
	\end{align*}
	see~\eqref{eq:Yk0}--\eqref{eq:Ykr}.

	In view of~\eqref{MomGnull1}, it remains to show that, as~$n\to \infty$,
	\begin{align}
		\label{tofinishproof31}
		{\rm E}_{\kappa_n,f}\n[
		g_{r,k}(\Ub) g_{s,\ell}(\Ub)
		]
		=
		{\rm E}_{0}\n[
		g_{r,k}(\Ub) g_{s,\ell}(\Ub)
		]
		+
		o(1)
	\end{align}
	for any~$k,\ell\in\{0,1,\ldots,\}$, $r\in\{1,\ldots,d_{p,k}\}$, and~$s\in\{1,\ldots,d_{p,\ell}\}$ (recall that~$d_{p,0}=1$ and~$g_{1,0}\equiv 1$). To do so, note that the continuity of~$f$ at zero implies that
	$$
	\sup_{\ub \in \mathcal{S}^{p-1}}
	|f(\kappa_n \ub\pr \thetab)-f(0)|
	=
	\sup_{s\in[-1,1]}
	|f(\kappa_n s)-f(0)|
	\to
	0
	,
	$$
	hence that
	\begin{align*}
		|c_{p,\kappa_n, f}^{-1}-c_p^{-1}|
		&\leq
		\int_{-1}^1 |f(\kappa_n s)-f(0)| (1-s^2)^{(p-3)/2}\,\mathrm{d}s
		\\
		&\leq
		\sup_{s\in[-1,1]}
		|f(\kappa_n s)-f(0)|
		\int_{-1}^1 (1-s^2)^{(p-3)/2}\,\mathrm{d}s
		\to
		0
	\end{align*}
	as~$n$ diverges to infinity. Therefore,
	$$
	\sup_{\ub \in \mathcal{S}^{p-1}}
	|c_{p,\kappa_n, f} f(\kappa_n \ub\pr \thetab)-c_p f(0)|
	\to
	0
	$$
	as~$n$ diverges to infinity.
	It then follows that
	\begin{align*}
		\bigg|
		\int_{\mathcal{S}^{p-1}} &
		g_{r,k}(\ub) g_{s,\ell}(\ub)
		(c_{p,\kappa_n, f} f(\kappa_n \ub\pr \thetab)- c_p f(0)) \,\mathrm{d}\sigma_p(\ub)
		\bigg|
		\nonumber
		\\
		&
		\leq
		\sup_{\ub \in \mathcal{S}^{p-1}}
		|c_{p,\kappa_n, f} f(\kappa_n \ub\pr \thetab)-c_p f(0)|
		\int_{\mathcal{S}^{p-1}}
		|g_{r,k}(\ub) g_{s,\ell}(\ub) |
		\,
		\mathrm{d}\sigma_p(\ub)
		\to 0
	\end{align*}
	as~$n$ diverges to infinity, which establishes~\eqref{tofinishproof31} and hence proves the result.
\end{proof}

The proof of Proposition~\ref{Momentprop} requires the following preliminary result.

\begin{lemma}
	\label{lemcontig}
	Fix integers~$p\geq 2$, $k>0$, and~$m\geq 0$. Let $g:\R\to \R$ be $k$ times differentiable at zero and let $(\kappa_n)$ be a positive real sequence that is~$o(1)$. Then, as~$n$ diverges to infinity,
	$$
	R_n^{(m)}(g)
	:=
	c_{p}
	\int_{-1}^1 (1-s^2)^{(p-3)/2} s^m g(\kappa_n s)  \,\mathrm{d}s
	=
	\sum_{\ell=0}^k \frac{a_{m+\ell}}{\ell !} \kappa_n^\ell g^{\underline{\ell}}(0)
	+o({\kappa_n^k}).
	$$
\end{lemma}

\begin{proof}{ of Lemma~\ref{lemcontig}}
	First note that \eqref{foref1} yields
	$$
	R_n^{(m)}(g)
	-
	\sum_{\ell=0}^{k-1} \frac{a_{m+\ell}}{\ell !} \kappa_n^\ell g^{\underline{\ell}}(0)
	=
	c_{p}
	\int_{-1}^1 (1-s^2)^{(p-3)/2}
	s^m
	\bigg\{
	g(\kappa_n s)
	- \sum_{\ell=0}^{k-1} \frac{(\kappa_n s)^\ell}{\ell !} g^{\underline{\ell}}(0)
	\bigg\}
	\,\mathrm{d}s
	.
	$$
	Letting~$t=\kappa_n s$ and using \eqref{foref1} again then provides
	$$
	R_n^{(m)}(g)
	-
	\sum_{\ell=0}^{k-1} \frac{a_{m+\ell}}{\ell !} \kappa_n^\ell g^{\underline{\ell}}(0)
	=
	{a_{m+k} \kappa_n^k}
	\int_{-\infty}^{\infty} h_n(t)
	\bigg( \frac{g(t)- \sum_{\ell=0}^{k-1}  \frac{t^\ell}{\ell!} g^{\underline{\ell}}(0)}{t^{k}} \bigg)
	\,\mathrm{d}t
	,
	$$
	where~$h_n$ is defined through
	$$
	t
	\mapsto
	h_n(t)
	=
	\frac{(\textstyle\frac{t}{\kappa_n})^{m+k}(1-(\textstyle\frac{t}{\kappa_n})^2)^{(p-3)/2}}{\int_{-\kappa_n}^{\kappa_n}(\textstyle\frac{s}{\kappa_n})^{m+k} (1-(\textstyle\frac{s}{\kappa_n})^2)^{(p-3)/2}\,\mathrm{d}s}
	\,
	\mathbb{I}[|t|\leq \kappa_n]
	.
	$$
	It can be checked that, since $\kappa_n=o(1)$, the~$h_n$'s form an \emph{approximate $\delta$-sequence}, in the sense that
	$
	\int_{-\infty}^\infty h_n(t)\,\mathrm{d}t=1
	$
	for all~$n$ and
	$
	\int_{-\varepsilon}^\varepsilon h_n(t)\,\mathrm{d}t \to 1
	$
	for any~$\varepsilon>0$. Hence,
	$$
	R_n^{(m)}(g)
	-
	\sum_{\ell=0}^{k-1} \frac{a_{m+\ell}}{\ell !} \kappa_n^\ell g^{\underline{\ell}}(0)
	=
	{a_{m+k} \kappa_n^k}
	\,
	\lim_{t\to 0} \bigg( \frac{g(t)-\sum_{\ell=0}^{k-1}  \frac{t^\ell}{\ell!} g^{\underline{\ell}}(0)}{t^k} \bigg)
	+
	o\big( {\kappa_n^k} \big)
	$$
	which yields the result by first using L'H\^{o}pital's rule~$k-1$ times and then the differentiability of~$g^{\underline{k-1}}$ at zero.
\end{proof}

The following result is a direct corollary of Lemma \ref{lemcontig} with~$m=0$.

\begin{lemma} \label{directLemma1}
	Fix integers~$p\geq 2$ and $k>0$. Let~$(\kappa_n)$ be a positive real sequence that is~$o(1)$ and~$f$ be an angular function that is~$k$ times differentiable at zero. Then, as~$n$ diverges to infinity,
	$$
	\frac{c_{p,\kappa_n, f}}{c_p}
	=
	\Bigg(
	\sum_{\ell=0}^k \frac{a_\ell}{\ell!}
	\kappa_n^\ell f^{\underline{\ell}}(0)
	+
	o({\kappa_n^k} )
	\Bigg)^{-1}
	,
	$$
	where~$c_{p,\kappa, f}$ is the normalizing constant in~\eqref{defcpkf}.
\end{lemma}

We can now prove Proposition \ref{Momentprop}.

\begin{proof}{ of Proposition \ref{Momentprop}}
	Writing
	$$
	e_{n,m}
	=
	\bigg(
	\frac{c_{p,\kappa_n, f}}{c_p}
	\bigg)
	c_{p}
	\int_{-1}^1  (1-s^2)^{(p-3)/2} s^m f(\kappa_n s) \,\mathrm{d}s
	$$
	and using Lemmas~\ref{lemcontig}--\ref{directLemma1} with~$k=q-m$ provides
	\begin{align*}
		e_{n,m}
		&=
		\Bigg(
		\sum_{\ell=0}^{q-m} \frac{a_\ell}{\ell!}
		\kappa_n^\ell f^{\underline{\ell}}(0)
		+
		o({\kappa_n^{q-m}})
		\Bigg)^{-1}
		\Bigg(
		\sum_{\ell=0}^{q-m} \frac{a_{m+\ell}}{\ell !} \kappa_n^\ell f^{\underline{\ell}}(0)
		+o({\kappa_n^{q-m}})
		\Bigg)
		\\
		&=
		\Bigg(
		\sum_{\ell=0}^{q-m} \beta_{\ell} \kappa_n^{\ell}
		+o({\kappa_n^{q-m}})
		\Bigg)^{-1}
		\Bigg(
		\sum_{\ell=0}^{q-m} \alpha_{\ell} \kappa_n^{\ell}
		+
		o({\kappa_n^{q-m}})
		\Bigg)
		,
	\end{align*}
	where we let
	$$
	\beta_{\ell}
	:=
	\frac{a_\ell}{\ell!}  f^{\underline{\ell}}(0)
	\quad
	\textrm{ and }
	\quad
	\alpha_{\ell}
	:=
	\frac{a_{\ell+m}}{\ell!}f^{\underline{\ell}}(0)
	.
	$$
	This yields
	$$
	e_{n,m}
	=
	\sum_{\ell=0}^{q-m}
	b_{m,\ell}
	\kappa_n^{\ell}
	+
	o({\kappa_n^{q-m}})
	,
	$$
	where the coefficients~$b_{m,\ell}$ satisfy
	$$
	\sum_{j=0}^\ell \beta_{\ell-j} b_{m,j}
	=
	\alpha_\ell
	,
	\quad \ell=0,1,\ldots,q-m
	,
	$$
	or equivalently, where~${\bf b}_m
	=
	(b_{m,0},b_{m,1},\ldots,b_{m,q-m})'
	$ satisfies~$\Ab_{q-m}{\bf b}_m=\vb_{q-m}^{(m)}$.
\end{proof}

The proof of Proposition~\ref{finalExprop} requires the following lemma.

\begin{lemma}
	\label{LemorthoGegen}
	For any positive integers~$j,k$ and for any integer~$p\geq 2$,
	$$
	c_p
	\int_{-1}^1
	C_{k}^{(p-2)/2}(s)
	C_{j}^{(p-2)/2}(s)
	(1-s^2)^{(p-3)/2}
	\,
	\mathrm{d}s
	=
	\frac{\delta_{k,j}}{t_{p,k}^2}
	,
	$$	
	where~$t_{p,k}$ was defined in Proposition~\ref{Exprop}. %
\end{lemma}

\begin{proof}{ of Lemma \ref{LemorthoGegen}}
	For~$p=2$, the result follows from the fact that~$c_2=1/\pi$ (see~\eqref{cpormula}) and the identity
	$$
	\int_{-1}^1
	C_{k}^{0}(s)
	C_{j}^{0}(s)
	(1-s^2)^{-1/2}
	\,
	\mathrm{d}s
	=
	\frac{\pi}{2}
	\delta_{j,k}
	;
	$$	
	see, e.g., Equation 3.11.8 in \cite{NIST:DLMF|SM}.
	For~$p=3$, Corollary~1.2.8 in \cite{Dai2013|SM} entails that
	\begin{align}
		c_p
		\int_{-1}^1
		C_{k}^{(p-2)/2}(s)
		&C_{j}^{(p-2)/2}(s)
		(1-s^2)^{(p-3)/2}
		\,
		\mathrm{d}s
		\nonumber
		\\
		& =
		\frac{c_p\omega_{p-1}}{\omega_{p-2}}
		\Big(1+ \frac{2k}{p-2}\Big)^{-1}
		C_{k}^{(p-2)/2}(1)
		\delta_{j,k}
		.
		\label{fks}
	\end{align}
	Using
	\begin{align*}
		g_{1,k}(\thetab_0)
		&=
		\frac{1}{\sqrt{d_{p,k}}}
		\Big(1+ \frac{2k}{p-2}\Big)
		C_k^{(p-2)/2}(1),
		\\
		g_{r,k}(\thetab_0)&=0
		,
		\quad
		r=2,\ldots,d_{p,k}
	\end{align*}
	in the addition formula~\eqref{iin} with~$\ub=\vb=\thetab_0$, we obtain the identity
	\begin{align}
		\label{Cinonerelation}
		C_k^{(p-2)/2}(1)
		=
		\Big(1+ \frac{2k}{p-2}\Big)^{-1}
		d_{p,k}
		.
	\end{align}
	Plugging this and the identity~$
	c_p\omega_{p-1}=\omega_{p-2}
	$ (which follows from~\eqref{cpormula}) into~\eqref{fks} then yields the result.
\end{proof}

We can now prove Proposition \ref{finalExprop}.

\begin{proof}{ of Proposition \ref{finalExprop}}
	Using~\eqref{Gegendefin}, Proposition~\ref{Momentprop} (with~$q=2(k-j)+r>k-2j=m$) yields
	\begin{align}
		{\rm E}\n_{\kappa_n,f}[C_k^{(p-2)/2}(\Ub_1\pr \thetab_0)]
		&=
		\sum_{j=0}^{\lfloor k/2 \rfloor} (-1)^{j} c^{(p-2)/2}_{k,j} e_{n,k-2j}
		\nonumber
		\\
		&=
		\sum_{j=0}^{\lfloor k/2 \rfloor} (-1)^{j} c^{(p-2)/2}_{k,j}
		\bigg(
		\sum_{\ell=0}^{k+r}
		b_{k-2j,\ell}
		\kappa_n^{\ell}
		+o\big({\kappa_n^{k+r}} \big)
		\bigg)
		\nonumber
		\\
		&=
		\sum_{\ell=0}^{k+r}
		w_\ell
		\kappa_n^{\ell}
		+o({\kappa_n^{k+r}})
		,
		\label{sss}
	\end{align}
	with
	$$
	w_\ell
	:=
	\sum_{j=0}^{\lfloor k/2 \rfloor}
	(-1)^{j} c^{(p-2)/2}_{k,j}
	b_{k-2j,\ell}
	.
	$$
	Letting~$\wb:=(w_0,w_1,\ldots,w_{k+r})'$, we have
	\begin{align}
		\wb
		&=
		\sum_{j=0}^{\lfloor k/2 \rfloor}
		(-1)^{j} c^{(p-2)/2}_{k,j}
		(b_{k-2j,0},b_{k-2j,1},\ldots,b_{k-2j,k+r})'
		\nonumber
		\\
		&=
		\Ab_{k+r}^{-1}
		\Bigg(
		\sum_{j=0}^{\lfloor k/2 \rfloor}
		(-1)^{j} c^{(p-2)/2}_{k,j}
		\vb_{k+r}^{(k-2j)}
		\Bigg)
		,
		\label{sss2}
	\end{align}
	where~$\Ab_{k+r}$ and~$\vb_{k+r}^{(k-2j)}$ are as in Proposition~\ref{Momentprop}. Now, for~$i=0,1,\ldots,k+r$,
	\begin{align*}
		\Bigg(
		\sum_{j=0}^{\lfloor k/2 \rfloor}
		(-1)^{j} c^{(p-2)/2}_{k,j}
		\vb_{k+r}^{(k-2j)}
		\Bigg)_{i}
		&=
		\sum_{j=0}^{\lfloor k/2 \rfloor}
		(-1)^{j}
		c^{(p-2)/2}_{k,j}
		\bigg(
		a_{k-2j+i} \frac{f^{\underline{i}}(0)}{i!}
		\bigg)
		\nonumber
		\\
		&
		=
		\frac{f^{\underline{i}}(0)}{i!}
		{\rm E}_0\n
		\Bigg[
		(\Ub_1\pr \thetab_0)^{i}
		\sum_{j=0}^{\lfloor k/2 \rfloor}
		(-1)^{j}
		c^{(p-2)/2}_{k,j}
		(\Ub_1\pr \thetab_0)^{k-2j}
		\Bigg]
		\nonumber
		\\
		&=
		\frac{f^{\underline{i}}(0)}{i!}
		{\rm E}_0\n\big[ (\Ub_1\pr \thetab_0)^{i} C_k^{(p-2)/2}(\Ub_1\pr \thetab_0)\big]
		.
	\end{align*}
	Recalling~\eqref{monomials}, this yields
	\begin{align*}
		\Bigg(
		\sum_{j=0}^{\lfloor k/2 \rfloor}
		(-1)^{j} c^{(p-2)/2}_{k,j}
		\vb_{k+r}^{(k-2j)}
		\Bigg)_{i}
		&
		=
		\frac{f^{\underline{i}}(0)}{i!}
		\sum_{j=0}^i m_{j,i}
		{\rm E}_0\n\big[ C_{j}^{(p-2)/2}(\Ub_1\pr \thetab_0) C_k^{(p-2)/2}(\Ub_1\pr \thetab_0)\big]
		\\
		&
		=
		\frac{f^{\underline{i}}(0)}{(i!)t_{p,k}^2}
		m_{k,i}
		\mathbb{I}[k\leq i]
		=
		\frac{1}{t_{p,k}^2}
		\Big({\bf z}_{k+r}^{(k)}\Big)_i
		,
	\end{align*}
	where we applied Lemma~\ref{LemorthoGegen}.
	Using this in~\eqref{sss}--\eqref{sss2}, we obtain that, with the vector~${\bf z}_{q-k}$ defined in the statement of the proposition,
	\begin{align}
		{\rm E}_{\kappa_n,f}\n\big[C_k^{(p-2)/2}(\Ub_1\pr \thetab_0)\big] &= \sum_{\ell=0}^{k+r}
		({\bf e}_{\ell+1, k+r+1}\pr {\bf w})
		\kappa_n^{\ell}
		+o({\kappa_n^{k+r}}) \nonumber \\
		&= \frac{1}{t_{p,k}^2}
		\sum_{\ell=k}^{k+r}
		({\bf e}_{\ell+1, k+r+1}\pr \Ab_{k+r}^{-1} {\bf z}_{k+r}^{(k)})
		\kappa_n^\ell+o({\kappa_n^{k+r}}), \label{anotherstep}
	\end{align}
	where the terms associated with~$\ell=0,1,\ldots,k-1$ cancel since~$\Ab_{k+r}^{-1}$ is a lower triangular matrix and the first~$k$ components of~${\bf z}_{k+r}^{(k)}$ are zero. Since~${\bf e}_{\ell+1,k+r+1}\pr \Ab_{k+r}^{-1}={\bf a}_{\ell+1;k+r}\pr$, the result is proved.
\end{proof}

%-------------------------------%
\section{Proofs of Section \ref{sec:localpowfiniteinf}}
%-------------------------------%

\begin{proof}{ of Proposition \ref{asymptnormtrunc}}
	For any~$i=1,\ldots,n$, consider the random vector~$\Zb_i:=(({\bf G}_{p,1}(\Ub_i))',\allowbreak \ldots, ({\bf G}_{p,K}(\Ub_i))')\pr$,  with dimension~$d:= d_{p,1}+\cdots+d_{p,K}$. Using the classical Cram\'{e}r--Wold device, we need to show that, for any $\wb \in \mathcal{S}^{d-1}$,
	$$
	\sqrt{n}\wb'({\bf T}_{p,K}\n- {\rm E}\n_{\kappa_n, f}[{\bf T}_{p,K}\n])
	=
	\frac{1}{\sqrt{n}}
	\sum_{i=1}^n
	(\wb\pr \Zb_i - {\rm E}\n_{\kappa_n, f}[\wb\pr \Zb_i ])
	=:
	\sum_{i=1}^n X_{ni}
	$$
	is asymptotically standard normal. From Proposition~\ref{Exprop}, we have that
	\begin{align*}
		s_n^2
		:=&\;
		{\rm Var}\n_{\kappa_n, f}\Bigg[\sum_{i=1}^n X_{ni}\Bigg]
		=
		n {\rm Var}\n_{\kappa_n, f}[X_{n1}]
		=
		\wb' {\rm Var}\n_{\kappa_n, f}[\Zb_1] \wb
		\\
		=&\;
		\wb' {\rm diag}({\bf I}_{d_{p,1}}, \ldots, {\bf I}_{d_{p,K}}) \wb + o(1)
		=1+o(1)
	\end{align*}
	as $n$ diverges to infinity. Now, since each component~$g_{r,k}$ of~${\bf G}_{p,k}$, $k=1,\ldots,K$, is a continuous function defined over the compact set~$\mathcal{S}^{p-1}$, there exists a positive constant~$C$ such that
	$$
	\max_{k=1,\ldots,K}
	\max_{r=1,\ldots,d_{p,k}}
	\sup_{u\in\mathcal{S}^{p-1}}
	|g_{r,k}(\ub)|
	\leq C
	,
	$$
	which implies that
	$$
	| X_{ni} |
	\leq
	\frac{1}{\sqrt{n}}
	\| \Zb_i - {\rm E}\n_{\kappa_n, f}[ \Zb_i ] \|
	\leq
	\frac{1}{\sqrt{n}}
	(\| \Zb_i \| + \| {\rm E}\n_{\kappa_n, f}[ \Zb_i ] \| )
	\leq
	\frac{2\sqrt{d}C}{\sqrt{n}}
	$$
	almost surely. Therefore,
	$$
	\frac{1}{s_n^3} \sum_{i=1}^n {\rm E}[| X_{ni}|^3]
	\leq
	\frac{8d^{3/2}C^3}{s_n^3\sqrt{n}}
	\to 0
	$$
	as $n$ diverges to infinity. The result thus follows from the Lyapounov central limit theorem; see, e.g., Theorem~5.2 in \cite{DasG2008|SM}.
\end{proof}

\begin{proof}{ of Theorem~\ref{Sobolfin}}
	By definition of~$k_v$ and~$K_v$,
	$$
	S^{(n)}_v
	=
	n
	({\bf T}_{p,k_v,K_v}\n)\pr
	{\rm diag}(v_{k_v}^2 {\bf I}_{d_{p,k_v}}, \ldots, v_{K_v}^2 {\bf I}_{d_{p,K_v}} )
	{\bf T}_{p,k_v,K_v}\n
	,
	$$
	with
	$$
	{\bf T}_{p,k,K}\n
	:=
	\frac{1}{n}
	\sum_{i=1}^n
	( ({\bf G}_{p,k}(\Ub_i))', \ldots, ({\bf G}_{p,K}(\Ub_i))' )\pr
	.
	$$
	Decompose then $\sqrt{n}{\bf T}_{p,k,K}\n$ into
	\begin{align*}
		\sqrt{n}{\bf T}_{p,k_v,K_v}\n
		&=
		\sqrt{n} \big( {\bf T}_{p,k_v,K_v}\n - {\rm E}\n_{\kappa_n,f}[{\bf T}_{p,k_v,K_v}\n] \big)
		+
		\sqrt{n} {\rm E}\n_{\kappa_n,f}[{\bf T}_{p,k_v,K_v}\n]\\
		&=:
		V_n+W_n
		,
	\end{align*}
	say. Now, it directly follows from Proposition~\ref{asymptnormtrunc} that~$V_n$ is asymptotically standard normal. For~$W_n$, we need to consider cases~\eqref{Sobolfin:1}--\eqref{Sobolfin:3} separately. In case~\eqref{Sobolfin:1}, Propositions~\ref{Exprop} and~\ref{finalExprop} entail that~$W_n=o(1)$ as~$\ny$, which establishes the result. In case~\eqref{Sobolfin:2}, under which~$\kappa_n = n^{-1/(2k_v)}\tau$, Propositions~\ref{Exprop} and~\ref{finalExprop} provide
	\begin{align}
		\label{ssTh41}
		W_n
		=
		\frac{m_{k_v,k_v} f^{\underline{k_v}}(0)}{(k_v!)t_{p,k_v}}
		\tau^{k_v}
		\binom{{\bf e}_{1,d_{p,k_v}}}{\bf 0}
		+
		o(1)
	\end{align}
	as $\ny$. The result then follows from the fact that, for any positive integer~$k$,
	$$
	m_{k,k}
	=
	\frac{1}{c^{0}_{k,0}}
	=
	\frac{1}{2^{k-1}}
	$$
	for $p=2$, and
	$$
	m_{k,k}
	=
	\frac{1}{c^{(p-2)/2}_{k,0}}
	=
	\frac{2^{-k} k!\Gamma(\frac{p-2}{2})}{\Gamma(\frac{p-2}{2}+k)}
	=
	\frac{k!}{2^{k} (\frac{p-2}{2})_{k}}
	$$
	for $p\geq 3$; see~\eqref{Gegendefin}, \eqref{Chebydefin}, and~\eqref{monomials}. Finally, in case~\eqref{Sobolfin:3},
	Propositions~\ref{Exprop} and~\ref{finalExprop} now entail
	\begin{align}
		\label{ssTh41iii}
		W_n
		=
		r_n
		\binom{{\bf e}_{1,d_{p,k_v}}}{\bf 0}
		,
	\end{align}
	where~$|r_n|\to\infty$ as $\ny$. This establishes the result.
\end{proof}

\begin{proof}{ of Theorem~\ref{Sobolfinbisbis}}
	\eqref{Sobolfinbisbis:1} For any~$k\in\{k_v,k_v+1,\ldots,K_v\}$, applying Proposition~\ref{Exprop}, then Proposition~\ref{finalExprop} with~$r=q-k$, yields
	\begin{align}
		\mathrm{E}\n_{\kappa_n,f} [{\bf G}_{p,k}(\Ub_1)]
		& =
		\frac{1}{t_{p,k}}
		\sum_{\ell=k}^{q}
		({\bf a}_{\ell+1;q}\pr {\bf z}^{(k)}_{q})
		\kappa_n^\ell
		{\bf e}_1
		+
		o(\kappa_n^{q})
		\nonumber
		\\
		& =
		\frac{1}{t_{p,k}}
		\sum_{\ell=k}^q
		\;
		\sum_{i=k}^\ell
		\frac{1}{i!} (\Ab_q^{-1})_{\ell+1,i+1} m_{k,i} f^{\underline{i}}(0)
		\kappa_n^\ell {\bf e}_{1}
		+
		o(\kappa_n^{q})
		,
		\label{dhks}
	\end{align}
	where the sum over~$i$ stops at~$\ell$ due to the lower triangular nature of~$\Ab_q^{-1}$. As noted below Proposition~\ref{finalExprop}, ${\bf A}_q^{-1}= \sum_{j=0}^{q} (-1)^j \Lb_q^j$, where~$\Lb_q$ is a strictly lower triangular matrix. Since~$(\Lb_q)_{rs}$ is zero if~$r-s$ is odd, it follows that a necessary condition for~$(\Ab_q^{-1})_{\ell+1,i+1}\neq 0$ is that~$\ell\sim i$. Recalling that, similarly, a necessary condition for~$m_{k,i}\neq 0$ is that~$k\sim i$,~\eqref{dhks} rewrites
	\begin{align}
		\mathrm{E}\n_{\kappa_n,f} [{\bf G}_{p,k}(\Ub_1)]
		=
		\frac{1}{t_{p,k}}
		\sum_{\substack{k \leq \ell \leq q \\ \ell \sim k}}
		s_{k,\ell}
		\kappa_n^\ell {\bf e}_{1}
		+
		o(\kappa_n^{q})
		,
		\label{keyexpr1}
	\end{align}
	with
	$$
	s_{k,\ell}
	=
	s_{k,\ell,q}
	:=
	\sum_{\substack{k \leq i \leq \ell \\ i \sim k}}
	\frac{1}{i!} (\Ab_q^{-1})_{\ell+1,i+1} m_{k,i} f^{\underline{i}}(0)
	.
	$$
	We then consider three cases.
	(A) Fix~$k\in\{k_v,\ldots,k_\dagger-1\}\,\cap\, {\rm Supp}(v)$, where~${\rm Supp}(v):=\{k:v_k\neq 0\}$. Then, there cannot exist~$i\in\{k,k+1,\ldots,k_*\}$ such that~$i\sim k$ and~$f^{\underline{i}}(0)\neq 0$ (the existence of such an integer~$i$ would indeed contradict the definition of~$k_*$ for~$i<k_*$ and the definition of~$k_\dagger$ for~$i=k_*$). It follows that
	$$
	\mathrm{E}\n_{\kappa_n,f} [{\bf G}_{p,k}(\Ub_1)]=O(\kappa_n^{k_*+1})
	,
	$$
	so that~$\sqrt{n}\mathrm{E}\n_{\kappa_n,f} [{\bf G}_{p,k}(\Ub_1)]=o(1)$.
	(B) Fix $k\in\{k_\dagger,k_\dagger+1,\ldots,k_*\}\cap {\rm Supp}(v)$. The same argument as in Case~(A) implies that there is no~$i\in\{k,k+1,\ldots,k_*-1\}$ such that~$i\sim k$ and~$f^{\underline{i}}(0)\neq 0$. Therefore,
	$$
	\mathrm{E}\n_{\kappa_n,f} [{\bf G}_{p,k}(\Ub_1)]
	=
	\frac{1}{t_{p,k}}
	s_{k,k_*}
	\kappa_n^{k_*}
	\mathbb{I}[k\sim k_*]
	{\bf e}_{1}
	+
	o(\kappa_n^{k_*})
	,
	$$
	where
	$
	s_{k,k_*}=(k_*!)^{-1}
	m_{k,k_*} f^{\underline{k_*}}(0)
	$ (recall that the diagonal entries of~$\Ab_q^{-1}$ are equal to one),
	which yields
	$$
	\sqrt{n}\mathrm{E}\n_{\kappa_n,f} [{\bf G}_{p,k}(\Ub_1)]
	=
	\frac{m_{k,k_*} f^{\underline{k_*}}(0) \tau^{k_*}}{(k_*!) t_{p,k}}
	\mathbb{I}[k\sim k_*]
	{\bf e}_{1}
	+
	o(1)
	.
	$$
	(C) Fix $k\in\{k_*+1,\ldots,K_v\}\, \cap\, {\rm Supp}(v)$. Then, \eqref{keyexpr1} trivially implies that $\mathrm{E}\n_{\kappa_n,f} [{\bf G}_{p,k}(\Ub_1)]=o(\kappa_n^{k_*})$, so that
	$$
	\sqrt{n}\mathrm{E}\n_{\kappa_n,f} [{\bf G}_{p,k}(\Ub_1)]
	=
	o(1)
	.
	$$
	Summing up,  for any~$k\in{\rm Supp}(v)$,
	$$
	\sqrt{n}\mathrm{E}\n_{\kappa_n,f} [{\bf G}_{p,k}(\Ub_1)]
	=
	\frac{m_{k,k_*} f^{\underline{k_*}}(0)\tau^{k_*}}{(k_*!)t_{p,k}}
	\mathbb{I}[k\sim k_*,k_\dagger\leq k\leq k_*]
	{\bf e}_{1}
	+
	o(1)
	.
	$$
	The weak limit result in Part~\eqref{Sobolfinbisbis:1} of the theorem then follows as in the proof of Theorem~\ref{Sobolfin}\eqref{Sobolfin:2}. Since Part~\eqref{Sobolfinbisbis:2} follows in the exact same way, it only remains to show that~$\xi_{p,k,k_*}(\tau)$ admits the alternative expression provided in the statement of the theorem. Clearly, it is sufficient to show that
	\begin{align}
		\label{rewritingxi}
		\frac{m_{k,k_*}}{t_{p,k}^2}
		=
		\sum_{j=0}^{\lfloor k/2 \rfloor}
		(-1)^{j}
		c^{(p-2)/2}_{k,j}
		a_{k+k_*-2j}
		.
	\end{align}
	To do so, note that, for any integer~$k_*>0$ and any~$k\in\{1,\ldots,k_*\}$, Lemma~\ref{LemorthoGegen} yields
	\begin{align*}
		c_p
		\int_{-1}^1
		&C_{k}^{(p-2)/2}(s)
		s^{k_*}
		(1-s^2)^{(p-3)/2}
		\,
		\mathrm{d}s
		\\
		& =
		c_p
		\sum_{\ell=0}^{k_*} m_{\ell,k_*}
		\int_{-1}^1
		C_{k}^{(p-2)/2}(s)
		C_{\ell}^{(p-2)/2}(s)
		(1-s^2)^{(p-3)/2}
		\,
		\mathrm{d}s
		\\
		& =
		\sum_{\ell=0}^{k_*}
		\frac{m_{\ell,k_*}\delta_{k,\ell}}{t_{p,k}^2}
		=
		\frac{m_{k,k_*}}{t_{p,k}^2}
		\cdot
	\end{align*}
	Now, using~\eqref{Gegendefin} and~\eqref{Chebydefin} (for~$p\geq 3$ and~$p=2$, respectively), then~\eqref{foref1}, we also have
	\begin{align*}
		c_p
		\int_{-1}^1&
		C_{k}^{(p-2)/2}(s)
		s^{k_*}
		(1-s^2)^{(p-3)/2}
		\,
		\mathrm{d}s
		\\
		& =
		\sum_{j=0}^{\lfloor k/2 \rfloor}
		(-1)^{j}
		c^{(p-2)/2}_{k,j}
		\bigg(
		c_p
		\int_{-1}^1
		s^{k-2j}
		s^{k_*}
		(1-s^2)^{(p-3)/2}
		\,
		\mathrm{d}s
		\bigg)
		\\
		& =
		\sum_{j=0}^{\lfloor k/2 \rfloor}
		(-1)^{j}
		c^{(p-2)/2}_{k,j}
		a_{k+k_*-2j}
		,
	\end{align*}
	which establishes~\eqref{rewritingxi}, hence the theorem.
\end{proof}

\begin{proof}{ of Theorem~\ref{Sobolfininfin}}
	We start with the proof of Part~\eqref{Sobolfininfin:2} of the result. For any integer~$K\geq k_v$, consider the truncated Sobolev test statistic
	$$
	S^{(n)}_{v,K}
	:=
	\sum_{k=k_v}^{K}
	v_k^2
	\bigg\| \frac{1}{\sqrt{n}} \sum_{i=1}^n {\bf G}_{p,k}(\Ub_i) \bigg\|^2
	;
	$$
	see~\eqref{Sobolgen} and~\eqref{untruncSobol}.
	We first show that there exists a positive integer~$N$ such that the Sobolev statistic~$S^{(n)}_v$ in~\eqref{untruncSobol} satisfies
	\begin{align}
		\label{unifconvK}
		\sup_{n\geq N}
		{\rm E}\n_{\kappa_n,f}[|S^{(n)}_v - S^{(n)}_{v,K}|]
		=
		o(1)
	\end{align}
	as~$K$ diverges to infinity. To do so, first note that
	\begin{align*}
		{\rm E}\n_{\kappa_n,f}[|S^{(n)}_v - S^{(n)}_{v,K}|]
		&=
		{\rm E}\n_{\kappa_n,f}
		\Bigg[
		\sum_{k=K+1}^{\infty}
		v_k^2
		\bigg\| \frac{1}{\sqrt{n}} \sum_{i=1}^n {\bf G}_{p,k}(\Ub_i) \bigg\|^2
		\Bigg]
		\\
		&=
		\sum_{k=K+1}^{\infty}
		v_k^2
		{\rm E}\n_{\kappa_n,f}
		\Bigg[
		\frac{1}{n}
		\sum_{i,j=1}^n
		({\bf G}_{p,k}(\Ub_i))' {\bf G}_{p,k}(\Ub_j)
		\Bigg]
		,
	\end{align*}
	where we used the monotone convergence theorem. Hence,
	\begin{align*}
		{\rm E}\n_{\kappa_n,f}&[|S^{(n)}_v - S^{(n)}_{v,K}|]\\
		&=
		\sum_{k=K+1}^{\infty} v_k^2
		\big\{
		{\rm E}\n_{\kappa_n,f}[\| {\bf G}_{p,k}(\Ub_1) \|^2]
		+
		(n-1)
		({\rm E}\n_{\kappa_n,f}[{\bf G}_{p,k}(\Ub_1)])\pr
		{\rm E}\n_{\kappa_n,f}[{\bf G}_{p,k}(\Ub_2)]
		\big\}
		\\
		&=
		\sum_{k=K+1}^{\infty} v_k^2
		\big\{
		d_{p,k}
		+
		(n-1)
		t_{p,k}^2
		\,
		\big(
		\mathrm{E}_{\kappa_n,f}\n[C_k^{(p-2)/2}(\Ub_1'\thetab_0)]
		\big)^2
		\big\}
		,
	\end{align*}
	where the last equality follows from~\eqref{iin}, \eqref{Cinonerelation}, and Proposition~\ref{Exprop}. Now, for~$k\geq k_v$, one has
	\begin{align*}
		\mathrm{E}_{\kappa_n,f}\n&[C_k^{(p-2)/2}(\Ub_1'\thetab_0)]
		\\
		& =
		\mathrm{E}_{\kappa_n,f}\n[C_k^{(p-2)/2}(\Ub_1'\thetab_0)]
		-
		\frac{c_{p,\kappa_n, f}}{c_p}
		\sum_{\ell=0}^{k_v-1}
		\frac{f^{\underline{\ell}}(0)}{\ell!}
		\kappa_n^\ell
		\mathrm{E}_{0}\n[ (\Ub_1'\thetab_0)^\ell C_k^{(p-2)/2}(\Ub_1'\thetab_0)]
		\\
		& =
		c_{p,\kappa_n, f}\int_{-1}^1
		C_{k}^{(p-2)/2}(s)
		(1-s^2)^{(p-3)/2}
		\bigg\{
		f(\kappa_n s)
		-
		\sum_{\ell=0}^{k_v-1}
		\frac{f^{\underline{\ell}}(0)}{\ell!} \kappa_n^\ell s^\ell
		\bigg\}
		\,
		\mathrm{d}s
		,
	\end{align*}
	which, by using the Cauchy--Schwartz inequality, yields
	\begin{align*}
		\big(\mathrm{E}_{\kappa_n,f}\n&[C_k^{(p-2)/2}(\Ub_1'\thetab_0)]\big)^2
		\\
		&
		\leq
		c_{p,\kappa_n, f}^2
		\bigg(
		\int_{-1}^1
		(C_{k}^{(p-2)/2}(s))^2
		(1-s^2)^{(p-3)/2}
		\,\mathrm{d}s
		\bigg)
		\\
		&
		\times
		\bigg(
		\int_{-1}^1
		(1-s^2)^{(p-3)/2}
		\bigg\{
		f(\kappa_n s)
		-
		\sum_{\ell=0}^{k_v-1}
		\frac{f^{\underline{\ell}}(0)}{\ell!} \kappa_n^\ell s^\ell
		\bigg\}^2
		\,
		\mathrm{d}s
		\bigg)
		.
	\end{align*}
	Lemma~\ref{LemorthoGegen} and then the fact that~$c_{p,\kappa_n, f}\to c_{p}$ as~$n\to\infty$ (see the proof of Proposition \ref{Exprop})
	ensure that, for~$n\geq N$ large enough (with~$N$ not depending on~$k$),
	\begin{align*}
		\big(\mathrm{E}_{\kappa_n,f}\n&[C_k^{(p-2)/2}(\Ub_1'\thetab_0)]\big)^2
		\\
		& \leq
		\frac{c_{p,\kappa_n, f}^2}{t_{p,k}^2c_{p}}
		\int_{-1}^1
		(1-s^2)^{(p-3)/2}
		\bigg(
		f(\kappa_n s)
		-
		\sum_{\ell=0}^{k_v-1}
		\frac{f^{\underline{\ell}}(0)}{\ell!} \kappa_n^\ell s^\ell
		\bigg)^2
		\,
		\mathrm{d}s
		\\
		& \leq
		\frac{2c_{p}}{t_{p,k}^2}
		\int_{-1}^1
		(1-s^2)^{(p-3)/2}
		\bigg(
		f(\kappa_n s)
		-
		\sum_{\ell=0}^{k_v-1}
		\frac{f^{\underline{\ell}}(0)}{\ell!} \kappa_n^\ell s^\ell
		\bigg)^2
		\,
		\mathrm{d}s
		.
	\end{align*}
	Proceeding as in the proof of Lemma~\ref{lemcontig}, we then obtain
	\begin{align*}
		(n-1)
		t_{p,k}^2&
		\big(\mathrm{E}_{\kappa_n,f}\n[C_k^{(p-2)/2}(\Ub_1'\thetab_0)]\big)^2
		\\
		&
		\leq
		2nc_{p}
		\int_{-1}^1
		(1-s^2)^{(p-3)/2}
		\bigg(
		f(\kappa_n s)
		-
		\sum_{\ell=0}^{k_v-1}
		\frac{f^{\underline{\ell}}(0)}{\ell!} \kappa_n^\ell s^\ell
		\bigg)^2
		\,
		\mathrm{d}s
		\\
		&
		\leq
		2n c_{p}a_{2k_v}
		\kappa_n^{2k_v}
		\int_{-1}^1
		h_n(t)
		\bigg(
		\frac{
			f(t)
			-
			\sum_{\ell=0}^{k_v-1}
			\frac{f^{\underline{\ell}}(0)}{\ell!} t^\ell
		}
		{
			t^{k_v}
		}
		\bigg)^2
		\,
		\mathrm{d}t
		,
	\end{align*}
	where the function~$h_n$, defined through
	$$
	t
	\mapsto
	h_n(t)
	=
	\frac{(\textstyle\frac{t}{\kappa_n})^{2k_v}(1-(\textstyle\frac{t}{\kappa_n})^2)^{(p-3)/2}}{\int_{-\kappa_n}^{\kappa_n}(\textstyle\frac{s}{\kappa_n})^{2k_v} (1-(\textstyle\frac{s}{\kappa_n})^2)^{(p-3)/2}\,\mathrm{d}s}
	\,
	\mathbb{I}[|t|\leq \kappa_n]
	,
	$$
	is still an approximate $\delta$-sequence in the sense described in the proof of Lemma~\ref{lemcontig}. Therefore, L'H\^opital's rule yields that
	\begin{align*}
		\int_{-1}^1
		h_n(t)
		\bigg(
		\frac{
			f(t)
			-
			\sum_{\ell=0}^{k_v-1}
			\frac{f^{\underline{\ell}}(0)}{\ell!} t^\ell
		}
		{
			t^{k_v}
		}
		\bigg)^2
		\,
		\mathrm{d}t
		\to
		\lim_{n\to\infty}
		\bigg(
		\frac{
			f(t)
			-
			\sum_{\ell=0}^{k_v-1}
			\frac{f^{\underline{\ell}}(0)}{\ell!} t^\ell
		}
		{
			t^{k_v}
		}
		\bigg)^2
		=
		\frac{(f^{\underline{k_v}}(0))^2}{(k_v!)^2}
	\end{align*}
	as~$n$ diverges to infinity. Thus, for any~$n\geq N$ with~$N$ large enough (still independent on~$k$), we have that
	\begin{align*}
		(n-1)
		t_{p,k}^2
		\big(\mathrm{E}_{\kappa_n,f}\n[C_k^{(p-2)/2}(\Ub_1'\thetab_0)]\big)^2
		&\leq
		\frac{3 c_{p} a_{2k_v} (f^{\underline{k_v}}(0))^2}{(k_v!)^2}
		n\kappa_n^{2k_v}
		\\
		&=
		\frac{3 c_{p} a_{2k_v} (f^{\underline{k_v}}(0))^2}{(k_v!)^2}
		\tau^{2k_v}
		.
	\end{align*}
	Consequently, for any~$n\geq N$,
	\begin{align*}
		{\rm E}\n_{\kappa_n,f}[|S^{(n)}_v - S^{(n)}_{v,K}|]
		&\leq
		\sum_{k=K+1}^{\infty} v_k^2
		\bigg(
		d_{p,k}
		+
		\frac{3 c_{p} a_{2k_v} (f^{\underline{k_v}}(0))^2}{(k_v!)^2}
		\tau^{2k_v}
		\bigg)
		\\
		&\leq
		\bigg(
		1
		+
		\frac{3 c_{p} a_{2k_v} (f^{\underline{k_v}}(0))^2}{(k_v!)^2}
		\tau^{2k_v}
		\bigg)
		\sum_{k=K+1}^{\infty} v_k^2
		d_{p,k}
		,
	\end{align*}
	which, in view of the summability condition assumed on the coefficients~$v_k$, establishes~\eqref{unifconvK}.
	
	Now, fix~$\varepsilon>0$ and pick~$K_0(\geq k_v)$ large enough to have
	$$
	\sup_{n\geq N}
	{\rm E}\n_{\kappa_n,f}[|S^{(n)}_v - S^{(n)}_{v,K_0}|]
	\leq \frac{\varepsilon}{3}
	\quad
	\textrm{ and }
	\quad
	\sum_{k=K_0+1}^{\infty}
	v_k^2
	d_{p,k}
	\leq \frac{\varepsilon}{3}
	$$
	(existence of~$K_0$ follows from~\eqref{unifconvK} and from the summability condition on the~$v_k$'s). Denoting as~$\mathcal{L}^{(n)}_{v}$ and~$\mathcal{L}^{(n)}_{v,K_0}$ the distributions, under~${\rm P}\n_{\kappa_n,f}$, of~$S^{(n)}_v$ and~$S^{(n)}_{v,K_0}$, respectively, we thus have that, for any~$n\geq N$,
	$$
	W_1(\mathcal{L}^{(n)}_{v},\mathcal{L}^{(n)}_{v,K_0})
	\leq
	{\rm E}\n_{\kappa_n,f}[|S^{(n)}_v - S^{(n)}_{v,K_0}|]
	\leq \frac{\varepsilon}{3}
	,
	$$
	where~$W_1(\mu,\nu)$ is the Wasserstein distance of order one between~$\mu$ and~$\nu$. Now, Theorem~\ref{Sobolfin}\eqref{Sobolfininfin:2} establishes that, as~$n\to\infty$,
	$$
	S^{(n)}_{v,K_0}
	\stackrel{\mathcal{D}}{\to}
	S^{\infty}_{v,K_0}
	:=
	\sum_{k=k_v}^{K_0} v_k^2 \chi^2_{d_{p,k}}(\delta_{k,k_v}\xi_{p,k_v}(\tau)).
	$$
	Since
	\begin{align*}
		{\rm E}\n_{\kappa_n,f}[|S^{(n)}_{v,K_0}|]
		&=
		{\rm E}\n_{\kappa_n,f}[S^{(n)}_{v,K_0}]
		=
		{\rm E}\n_{\kappa_n,f}
		\bigg[
		\sum_{k=k_v}^{K_0}
		v_k^2
		\bigg\| \frac{1}{\sqrt{n}} \sum_{i=1}^n {\bf G}_{p,k}(\Ub_i) \bigg\|^2
		\bigg]
		\\
		&=
		\sum_{k=k_v}^{K_0}
		v_k^2
		{\rm E}\n_{\kappa_n,f}
		\bigg[
		\frac{1}{n}
		\sum_{i,j=1}^n
		({\bf G}_{p,k}(\Ub_i))' {\bf G}_{p,k}(\Ub_j)
		\bigg]
		\\
		&=
		\sum_{k=k_v}^{K_0}
		v_k^2
		\big\{
		d_{p,k}
		+
		(n-1)
		({\rm E}\n_{\kappa_n,f}[{\bf G}_{p,k}(\Ub_1)])\pr
		{\rm E}\n_{\kappa_n,f}[{\bf G}_{p,k}(\Ub_2)]
		\big\}
		,
	\end{align*}
	using~\eqref{ssTh41} yields
	\begin{align*}
		{\rm E}\n_{\kappa_n,f}[|S^{(n)}_{v,K_0}|]
		&=
		\sum_{k=k_v}^{K_0}
		v_k^2
		\bigg\{
		d_{p,k}
		+
		\frac{m^2_{k_v,k_v}(f^{\underline{k_v}}(0))^2\tau^{2k_v}}{(k_v!)^2 t^2_{p,k_v}}
		\delta_{k,k_v}
		+
		o(1)
		\bigg\}
		\\
		&\to
		\sum_{k=k_v}^{K_0}
		v_k^2
		(
		d_{p,k}
		+
		\delta_{k,k_v}\xi_{p,k_v}(\tau)
		)
		=
		{\rm E}[|S^{\infty}_{v,K_0}|]
		,
	\end{align*}
	so that, in the sense of Definition~6.8 from \cite{Vil2009|SM}, $\mathcal{L}^{(n)}_{v,K_0}$ (the distribution of~$S^{(n)}_{v,K_0}$ under~${\rm P}\n_{\kappa_n,f}$) converges weakly to the distribution, $\mathcal{L}^{\infty}_{v,K_0}$, of~$S^{\infty}_{v,K_0}$ in the Wasserstein space of order one. Therefore, Theorem~6.9 from the same monograph implies that
	$$
	W_1(\mathcal{L}^{(n)}_{v,K_0},\mathcal{L}^{\infty}_{v,K_0})
	\leq \frac{\varepsilon}{3}
	$$
	for~$n$ large enough. Finally, denoting as~$\mathcal{L}^{\infty}_{v}$ the distribution of
	$$
	S^\infty_v
	:=
	\sum_{k=k_v}^\infty v_k^2 \chi^2_{d_{p,k}}(\delta_{k,k_v}\xi_{p,k_v}(\tau))
	,
	$$
	we have that
	$$
	W_1(\mathcal{L}^{\infty}_{v,K_0},\mathcal{L}^{\infty}_{v})
	\leq
	{\rm E}[| S^{\infty}_{v,K_0}-S^\infty_v |]
	=
	{\rm E}\Bigg[\sum_{k=K_0+1}^\infty v_k^2 \chi^2_{d_{p,k}}(\delta_{k,k_v}\xi_{p,k_v}(\tau))\Bigg]
	=
	\sum_{k=K_0+1}^\infty v_k^2
	d_{p,k}
	\leq
	\frac{\varepsilon}{3}
	\cdot
	$$
	Therefore, the triangle inequality yields that
	$$
	W_1(\mathcal{L}^{(n)}_v,\mathcal{L}^{\infty}_{v})
	\leq \varepsilon,
	$$
	for~$n$ large enough, which, using again Theorem~6.9 in \cite{Vil2009|SM}, establishes that~$S^{(n)}_v$ converges to~$S^\infty_v$ in distribution under~${\rm P}\n_{\kappa_n,f}$. This proves Part~\eqref{Sobolfininfin:2} of the result. The proof of Part~\eqref{Sobolfininfin:1} of the result follows entirely similarly from Theorem~\ref{Sobolfin}\eqref{Sobolfin:1} (the proof is actually slightly simpler since the expectation shift at rank~$k_v$ vanishes in the corresponding setup). Finally, note that, with the sequence~$\tilde{v}=(\tilde{v}_k):=(v_{k_v} \delta_{k,k_v})$, one trivially has that
	$$
	{\rm P}\n_{\kappa_n,f}[S^{(n)}_v>c]
	\geq
	{\rm P}\n_{\kappa_n,f}[S^{(n)}_{\tilde{v}}>c]
	$$
	for any~$c>0$ and any integer~$n$, so that Part~\eqref{Sobolfininfin:3} of the result follows from Theorem~\ref{Sobolfin}\eqref{Sobolfin:3}.
\end{proof}

%\bibliographystyle{apalike}
%\bibliography{DirectionalStats.bib,no-DirectionalStats.bib}
% DirectionalStats.bib v1.25 (2021-07-20)

\fi

\end{document}

%% file: preamble_arXiv.tex
\usepackage[utf8]{inputenc}
\usepackage[T1]{fontenc}

\usepackage{amsthm}
\usepackage{amsmath}
\usepackage{amsfonts}
\usepackage{amssymb}

\usepackage{graphicx} 
\usepackage{float}
\usepackage{booktabs}
\usepackage{multirow}
\usepackage{rotating}
\usepackage{array}
\usepackage{xcolor}
\usepackage{subcaption}
\usepackage[ruled]{algorithm2e}

\newcolumntype{R}[1]{>{\raggedleft\let\newline\\\arraybackslash\hspace{0pt}}m{#1}}
\newcolumntype{L}[1]{>{\raggedright\let\newline\\\arraybackslash\hspace{0pt}}m{#1}}

\usepackage{breakcites}

\usepackage{nowidow}

\usepackage{enumitem}

\usepackage[top=1.5cm,bottom=2cm,right=2.5cm,left=2.5cm]{geometry}
\usepackage{titling}

\usepackage{natbib}

\usepackage{ifplatform}

\usepackage{textcomp}

\ifwindows
\fi

\allowdisplaybreaks

\newcommand{\R}{\mathbb{R}}

\newcommand{\N}{\mathbb{N}}
\newcommand{\Sb}{\mathbf{S}}
\newcommand{\thetab}{\boldsymbol\theta}
\newcommand{\etab}{\boldsymbol\eta}
\newcommand{\xib}{\boldsymbol\xi}

\newcommand{\Zb}{\mathbf{Z}}

\newcommand{\ub}{\mathbf{u}}

\newcommand{\vb}{\mathbf{v}}
\newcommand{\xb}{\mathbf{x}}

\newcommand{\wb}{\mathbf{w}}
\newcommand{\Ab}{\mathbf{A}}

\newcommand{\Ob}{\mathbf{O}}

\newcommand{\Ub}{\mathbf{U}}

\newcommand{\Lb}{\mathbf{L}}

\DeclareFontFamily{OT1}{pzc}{}
\DeclareFontShape{OT1}{pzc}{m}{it}{<-> s * [1.10] pzcmi7t}{}
\DeclareMathAlphabet{\mathpzc}{OT1}{pzc}{m}{it}

\newtheorem{theorem}{Theorem}
\newtheorem{corollary}{Corollary}

\newtheorem{proposition}{Proposition}
\newtheorem{lemma}{Lemma}

\newtheorem{thm}{Theorem}
\renewenvironment{proof}[1]{\textit{Proof#1.}}{\qed\\} 

\newcommand{\n}{^{(n)}}
\newcommand{\ny}{n\rightarrow\infty}

\newcommand{\Gamb}{{\boldsymbol \Gamma}}